\documentclass[preprint,10pt]{article}
\usepackage[T1]{fontenc}
\usepackage{amsmath,amssymb}
\usepackage{mathptmx}
\usepackage{graphicx,graphics}
\usepackage{tikz}
\usepackage{hyperref}
\usepackage{amsthm}
\usepackage{floatrow}
\usepackage{subfig}
\newtheorem{theorem}{Theorem}
\newtheorem{corollary}{Corollary}[theorem]
\newtheorem{lemma}[theorem]{Lemma}
\theoremstyle{definition}
\newtheorem{definition}{Definition}[section]

\newtheorem{remark}{Remark}[section]
\newtheorem{prop}[theorem]{Proposition}
\usepackage{mathrsfs}
\usepackage{enumerate}
\usepackage{float}

\usepackage{tkz-graph}
\usepackage{subfig}

\usetikzlibrary{shapes.geometric,arrows}
\usetikzlibrary{mindmap}
\usetikzlibrary{arrows,decorations.markings}
\definecolor{qqccqq}{rgb}{0,0.8,0}
\definecolor{qqzzcc}{rgb}{0,0.6,0.8}
\tikzset{every loop/.style={min distance=10mm,looseness=10}}
\tikzset{every state/.style={minimum size=2mm}}

\newcommand\FCCO{\mathcal F_\mathrm{CCO}}

\newcommand\MI[1]{M_{\mathit{I}}(#1)}
\newcommand\MIV{M_{\mathit{IV}}}
\newcommand\MV{M_{\mathit{V}}}
\newcommand\miop{\mathbin{\oplus}}
\newcommand\MIast[1]{M_{\mathit{I}}^*(#1)}

\newcommand\MVast{M_{\mathit{V}}^*}
\newcommand\gc{\mathcal{G}_{\mathit{co}}^p}
\newcommand\gs{\mathcal{G}_{\mathit{co}}}
\newcommand\ForbRow{\mathcal F_{\textup{circR}}}
\newcommand\trans{^{\mathit t}}

\newcommand\Ztwo{\ensuremath{Z_2^*}}
\newcommand\Ztthree{\ensuremath{Z_3^*}}
\newcommand\Zfour{\ensuremath{Z_4^*}}
\newcommand\Zfive{\ensuremath{Z_5}}
\newcommand\ZfiveTrans{\ensuremath{Z_5\trans}}

\newcommand\CoZtwo{\ensuremath{\overline{Z_2^*}}}
\newcommand\CoZfour{\ensuremath{\overline{Z_4^*}}}

\numberwithin{equation}{section}
\title{Forbidden Induced Subgraph Characterization of Word-Representable Co-bipartite Graphs}
\author{
	Eshwar Srinivasan and Ramesh Hariharasubramanian \footnote{Corresponding Author}\\
	{\small Department of Mathematics,
	Indian Institute of Technology Guwahati,
	Assam - 781039, India}\\
	{\small s.eshwar@iitg.ac.in and ramesh\_h@iitg.ac.in}
}
\date{}
\begin{document}
	    \maketitle      
	\begin{abstract}
		A graph $G$ with vertex set $V(G)$ and edge set $E(G)$ is said to be \emph{word-representable} if there exists a word \(w\) over the alphabet \(V(G)\) such that, for any two distinct letters \(x,y \in V(G)\), the letters \(x\) and \(y\) alternate in \(w\) if and only if \((x,y) \in E(G)\). Equivalently, a graph is word-representable if and only if it admits a \emph{semi-transitive orientation}, that is, an acyclic orientation in which, for every directed path \(v_0 \rightarrow v_1 \rightarrow \cdots \rightarrow v_m\) with \(m \ge 2\), either there is no arc between \(v_0\) and \(v_m\), or, for all \(1 \le i < j \le m\), there exists an arc from \(v_i\) to \(v_j\). In this work, we provide a comprehensive structural and algorithmic characterization of word-representable \emph{co-bipartite graphs}, a class of graphs whose vertex set can be partitioned into two cliques.
		
		This work unifies graph-theoretic and matrix-theoretic perspectives. We first establish that a co-bipartite graph is a circle graph if and only if it is a permutation graph, thereby deriving a minimal forbidden induced subgraph characterization for co-bipartite circle graphs. The central contribution then connects semi-transitivity with the circularly compatible ones property of binary matrices. In addition to the structural characterization, the paper introduces a linear-time recognition algorithm for semi-transitive co-bipartite graphs, utilizing Safe’s matrix recognition framework.
		
	\end{abstract}

	\section{Introduction}
	
	A graph $G$ with vertex set $V(G)$ and edge set $E(G)$ is said to be \emph{word-representable} if there exists a word \(w\) over the alphabet \(V(G)\) such that for any two distinct letters \(x,y \in V\), the letters \(x\) and \(y\) alternate in \(w\) if and only if \((x,y)\) is an edge in \(E(G)\). This concept, which first emerged in the study of the Perkins semigroup through the work of Kitaev and Seif~\cite{kitaev2008word}, has generated significant interest and inspired a broad range of research focused on the recognition and characterization of such graphs.
	
	The fundamental structural characterization of word-representable graphs was given by Halld\'orsson, Kitaev, and Pyatkin~\cite{halldorsson2011alternation}, who showed that a graph is word-representable precisely when it admits a \emph{semi-transitive orientation}. A \textit{semi-transitive} orientation is an acyclic orientation such that, for any directed path $v_0 \rightarrow v_1 \rightarrow \cdots \rightarrow v_m$ with $m \ge 2$, either there is no arc between $v_0$ and $v_m$, or, for all $1 \le i < j \le m$, an arc exists from $v_i$ to $v_j$. This remarkable result shifted the study of word-representable graphs from a purely combinatorial setting—concerned with the arrangement of letters in words—to a graph-theoretic framework focused on acyclic orientations that avoid specific forbidden patterns known as \emph{shortcuts}. The semi-transitivity condition can be viewed as a natural relaxation of the transitivity property that defines comparability graphs, thereby revealing a rich hierarchy of graph classes connected through their orientability properties. For further reading on word-representable graphs, we refer the reader to~\cite{kitaev2008word, kitaev2008representable, kitaev2011representability, kitaev13, HALLDORSSON2016164, kitaev2015words, kitaev2017comprehensive, broere2018word, kitaev2021human, srinivasan2024minimum, huang2024embedding}.
	
	A \emph{hereditary graph class} is one that is closed under taking induced subgraphs and can often be described in terms of a family of \emph{forbidden induced subgraphs}. In other words, a graph belongs to the class if and only if it contains no member of a specified family \(\mathcal{F}\) as an induced subgraph. The family \(\mathcal{F}\) is referred to as the collection of \emph{minimal forbidden induced subgraphs} for the class. For instance, \emph{cographs} are exactly those graphs that contain no induced path on four vertices (\(P_4\))~\cite{corneil1981}, while \emph{trivially perfect graphs} are characterized by the absence of both the path \(P_4\) and the cycle \(C_4\)~\cite{golumbic1980}.
	
	The study of hereditary graph classes and their forbidden subgraph characterizations plays a central role in both structural graph theory and algorithmic design. When the collection \(\mathcal{F}\) of minimal forbidden induced subgraphs is finite, testing membership in the class amounts to verifying the absence of finitely many patterns, often leading to efficient, polynomial-time recognition algorithms. However, not every hereditary graph class admits such a finite forbidden subgraph characterization, which presents additional challenges for structural understanding and algorithmic treatment. 
	
	The class of word-representable graphs is hereditary, meaning that every induced subgraph of a word-representable graph is also word-representable. Because of this property, an important open problem is to identify all minimal forbidden induced subgraphs for this class. Although word-representable graphs can be characterized using semi-transitive orientations, a complete description in terms of forbidden induced subgraphs is still unknown. 
	
	A natural way to make progress on this problem is to study well-structured hereditary subclasses. Within these restricted settings, we can investigate whether word-representability admits a forbidden induced subgraph characterization and whether there exists an efficient algorithm to recognize such graphs.

	Let $\mathcal{E}_{i,j}$ denote the class of graphs whose vertex set can be partitioned into at most $i$ independent sets and $j$ cliques. For instance, $\mathcal{E}_{2,0}$ is exactly the class of bipartite graphs, while $\mathcal{E}_{3,0}$ coincides with the class of $3$-colorable graphs. Within the framework of graph classes, several boundary cases are well understood.  It is known that $\mathcal{E}_{i,j}$ is word-representable only for $(i,j) \in \{(1,0), (0,1), (2,0), (3,0)\}$~\cite{kitaev2015words}. This naturally leads to the question of determining the word-representability of $\mathcal{E}_{i,j}$ for other values of $i$ and $j$.
	
	In particular, for $\mathcal{E}_{1,1}$, which corresponds to the class of split graphs, a necessary and sufficient condition for word-representability was established by Kitaev \emph{et al.}~\cite{kitaev2021word}. Since then, significant attention has been given to the structural and algorithmic aspects of word-representable split graphs~\cite{chen2022representing, iamthong2021semi, iamthong2022word, kitaev2021word, kitaev2024semi, roy2025word}. Recently, we obtained a forbidden induced subgraph characterization for word-representable split graphs in~\cite{srinivasan2025forbidden}.
	
	Co-bipartite graphs ($\mathcal{E}_{0,2}$) are graphs whose vertex set can be divided into two cliques. They form the next dense hereditary class after split graphs. In bipartite graphs, the structure is sparse, which makes orientation constraints easier to handle. In contrast, co-bipartite graphs contain large cliques, so semi-transitive orientations must be chosen carefully to avoid creating shortcuts.
	
	Another important feature makes co-bipartite graphs especially interesting: their vertex partition naturally gives rise to a $(0,1)$-matrix representation. This matrix viewpoint allows us to translate orientation constraints into matrix properties. In particular, we show that semi-transitivity in co-bipartite graphs is exactly characterized by the circularly compatible ones property. Therefore, co-bipartite graphs offer a unified setting where orientation theory, permutation and circle graph structure, and circular-ones matrix theory come together in a natural and structured way.
	
	The word-representability of co-bipartite graphs was first investigated by Chen \emph{et al.}~\cite{chen2025word}, who characterized word-representable co-bipartite graphs in terms of forbidden induced subgraphs under the condition that one of the cliques has size at most four.
	
	Later, Das and Hariharasubramanian~\cite{das2025word} provided a complete characterization through necessary and sufficient conditions. In addition, they established an upper bound on the representation number of word-representable co-bipartite graphs~\cite{das2025representation}.

	Building on this growing body of research, the present work concentrates on a detailed structural and algorithmic characterization of word-representable co-bipartite graphs. In particular, we provide a \emph{forbidden induced subgraph characterization} for this class of graphs, grounded in the theory of circularly compatible ones matrices. This result unifies two previously separate lines of research: the geometric theory of circle graphs, developed notably by Bouchet~\cite{MR1256586}, and the matrix-theoretic framework for circular-ones properties, introduced by Tucker~\cite{MR2618994, MR309810} and later refined by Safe~\cite{mdsafe1, mdsafe2}.
	
	The central theorem of this paper establishes that a co-bipartite graph is semi-transitive if and only if its associated $(0,1)$-matrix satisfies the \emph{circularly compatible ones property}. This equivalence provides not only a minimal forbidden induced subgraph characterization but also a linear-time recognition algorithm for semi-transitive co-bipartite graphs.
	
	Using an algebraic characterization of circle graphs, Bouchet~\cite{MR1256586} proved that if a bipartite graph \(G\) is the complement of a circle graph, then \(G\) itself is also a circle graph. Later, Esperet and Stehl\'ik~\cite{MR4062292} provided an elementary proof of this result. In Section~\ref{s3}, we present an even more elementary and concise proof of a stronger version of this theorem: if a bipartite graph \(G\) is the complement of a circle graph, then both \(G\) and \(\overline{G}\) are permutation graphs. This result can also be viewed as a direct consequence of the proof by Esperet and Stehl\'ik~\cite{MR4062292}.
	
	The remainder of the paper is organized as follows. Section~\ref{s2} introduces the necessary definitions and preliminaries, including key concepts from graph theory, word-representability, and matrix properties. In Section~\ref{s3}, we present a forbidden induced subgraph characterization of co-bipartite circle graphs and establish their equivalence with co-bipartite permutation graphs. Section~\ref{s4} contains the main result, proving the equivalence between semi-transitivity and the circularly compatible ones property; by combining this equivalence with Safe’s forbidden submatrix characterization, we obtain both a forbidden induced subgraph characterization and a linear-time recognition algorithm. Finally, Section~\ref{s5} concludes the paper and outlines directions for future research.
	
	\section{Definitions and Preliminaries}\label{s2}
	In this section, we provide some definitions and preliminaries that will be useful in the remainder of the paper. For each positive integer $k$, we denote the set $\{1, 2, \ldots, k\}$ by $[k]$. If $k = 0$, we define $[k]$ to be the empty set.
	\subsection{Graph}
	Let $G$ be a graph with vertex set $V(G)$ and edge set $E(G)$. For a vertex $v \in V(G)$, the \emph{neighborhood} of $v$ in $G$, denoted by $N_G(v)$, is the set of vertices adjacent to $v$ in $G$. The \emph{closed neighborhood} of $v$, denoted by $N_G[v]$, is defined as $N_G(v) \cup \{v\}$. For any subset $X \subseteq V(G)$, the \emph{subgraph of $G$ induced by $X$} is the graph with vertex set $X$ and edge set consisting of all edges of $G$ whose endpoints both lie in $X$. If $H$ is a graph, we say that $G$ \emph{contains $H$ as an induced subgraph} if $H$ is isomorphic to some induced subgraph of $G$. 
	
	An \emph{independent set} (respectively, \emph{clique}) of a graph $G$ is a set of vertices that are pairwise nonadjacent (respectively, pairwise adjacent). A graph class is said to be \emph{hereditary} if it is closed under taking induced subgraphs. The \emph{complement} of a graph $G$, denoted by $\overline{G}$, is the graph with vertex set $V(G)$ and edge set $E(\overline{G}) = \{xy \colon xy \notin E(G)\}$.

	The \emph{union} of two graphs $G_1$ and $G_2$, denoted by $G_1 \cup G_2$, is the graph with vertex set $V(G_1) \cup V(G_2)$ and edge set $E(G_1) \cup E(G_2)$. The \emph{join} of $G_1$ and $G_2$, denoted by $G_1 + G_2$, is obtained from $G_1 \cup G_2$ by adding all edges between every vertex of $G_1$ and every vertex of $G_2$. A vertex $v$ of a graph $G$ is called \emph{universal} if it is adjacent to every other vertex of $G$, and \emph{isolated} if it is adjacent to none.
	
	A graph is called an \emph{intersection graph} if there exists a set of objects such that each vertex corresponds to an object, and two vertices are adjacent if and only if the corresponding objects have a non-empty intersection. A graph $G$ is said to be \emph{bipartite} if its vertex set can be partitioned into two subsets such that each subset induces an independent set. A graph is \emph{co-bipartite} if its complement is bipartite. A graph is called a \emph{circle graph} if its vertices can be associated with the chords of a circle such that two vertices are adjacent if and only if the corresponding chords intersect.  
	
	Let $\pi$ be a permutation over $[n]$. The \emph{permutation graph} of $\pi$ is the graph with vertex set $[n]$, where two vertices $i$ and $j$ with $i < j$ are adjacent if and only if $j$ occurs before $i$ in $\pi$. The class of permutation graphs forms a subclass of the class of circle graphs.
	
	A \emph{local complementation} of a graph $G$ at a vertex $v$ is the graph, denoted by $G \ast v$, obtained by replacing the subgraph induced by the neighbours of $v$ with its complement. A sequence of local complementations at vertices $v_1, v_2, \ldots, v_k$ is denoted by $G \ast v_1 v_2 \ldots v_k$. Two graphs $G$ and $G'$ are said to be \emph{locally equivalent} if there exists a sequence $m$ of vertices such that $G' = G \ast m$. An \emph{$\ell$-reduction} of $G$ is an induced subgraph of a graph locally equivalent to $G$.
	
	All graphs in this study are assumed to be simple, that is, finite, undirected, and having no loops or multiple edges. For additional graph-theoretic terminology and notation, we refer the reader to \cite{MR1367739}.
	\subsection{Sequences}
	Let $a = a_1 a_2 \ldots a_k$ be a sequence of length $k$. The \emph{shift} of $a$ is defined as the sequence $a_2 a_3 \ldots a_k a_1$. The length of any sequence $a$ is denoted by $|a|$. A \emph{binary bracelet}~\cite{MR1857399} is the lexicographically smallest sequence in an equivalence class of binary sequences under the operations of shifting and reversal. For each $k \geq 4$, let $A_k$ denote the set of binary bracelets of length $k$. In particular, let $A_3 = {000, 111}$. Although $001$ and $011$ are also binary bracelets of length $3$, they are not included in $A_3$.
	\subsection{Matrices}
	
	All matrices considered in this article are \emph{binary}, that is, every entry of each matrix is either $0$ or $1$. Each row $r$ of a matrix $M$ is defined as the set of columns of $M$ having a $1$ in row $r$. Let $M$ be an $m \times n$ matrix. The rows and columns of $M$ are labeled by the sets $\{r_1, r_2, \ldots, r_m\}$ and $\{c_1, c_2, \ldots, c_n\}$, respectively. The \emph{complement} of a row $r_i$ is obtained by interchanging $0$s and $1$s in that row. The \emph{complement} of $M$, denoted by $\overline{M}$, is the matrix obtained from $M$ by complementing every row of $M$. A row or a column is said to be \emph{trivial} if it consists entirely of $0$s or entirely of $1$s.
	
	If $a$ is a binary sequence of length $k$, we denote by $a \miop M$ the matrix obtained from $M$ by complementing all rows $r_i$, $i \in [k]$, such that $a_i = 1$. Let $M$ and $M'$ be matrices. The matrix $M$ is said to \emph{contain} $M'$ as a \emph{configuration} if some submatrix of $M$ equals $M'$ up to permutations of rows and columns. We say that $M$ and $M'$ \emph{represent the same configuration} if they are equal up to permutations of rows and columns. If $M$ is a matrix, then $M^*$ denotes the matrix obtained from $M$ by adding an empty column.
	
	\subsection{Algorithms}
	When a matrix $M$ is given, we define its \emph{size}, denoted by $\mathrm{size}(M)$, as the sum of its number of rows, its number of columns, and its number of $1$-entries. An algorithm whose input is a matrix $M$ is said to run in linear time if its running time is $O(\mathrm{size}(M))$. For algorithms that take a graph as input, we let $n$ and $m$ denote the numbers of vertices and edges, respectively. Such an algorithm is called linear-time if it runs in $O(n+m)$ time. We represent input matrices as lists of rows, where each row is stored as the list of column indices containing a $1$ in that row. Graphs are represented using adjacency lists. Under these conventions, a matrix and a graph require $O(\mathrm{size}(M))$ and $O(n+m)$ space, respectively.
	
	\subsection{Word-representability}
	
	Let $w$ be a word over some finite alphabet, and let $x$ and $y$ be two distinct letters appearing in $w$. We say that $x$ and $y$ \emph{alternate} in $w$ if, after deleting all letters of $w$ except $x$ and $y$, the resulting word is either of the form $xyxy\ldots$ (of odd or even length) or $yxyx\ldots$ (of odd or even length). Hence, by definition, if $w$ contains only one occurrence each of $x$ and $y$, then $x$ and $y$ alternate in $w$. The concept of alternation forms the combinatorial foundation of the notion of word-representability in graphs.
	
	\begin{definition}\label{def:wordrepresentable}
		A graph $G$ is said to be \emph{word-representable} if there exists a word $w$ over the alphabet $V(G)$ such that, for every pair of distinct vertices $x, y \in V(G)$, the edge $xy$ belongs to $E(G)$ if and only if $x$ and $y$ alternate in $w$. Such a word $w$ is called a \emph{word-representant} of $G$, and we say that $w$ \emph{represents} $G$. Furthermore, $w$ must contain every letter of $V(G)$ at least once.
	\end{definition}
	
	The class of word-representable graphs was introduced by Kitaev and his collaborators as a natural generalization of several well-studied graph classes, including circle graphs and comparability graphs. Word-representable graphs have deep connections to poset theory, graph orientations, and permutation patterns.
	
	The following notion of \emph{semi-transitive orientation}, introduced in~\cite{halldorsson2011alternation}, provides a structural characterization of word-representable graphs.
	
	\begin{definition}\label{def:st_orientation}
		An acyclic orientation of a graph $G$ is said to be \emph{semi-transitive} if, for every directed path
		$v_0 \rightarrow v_1 \rightarrow \cdots \rightarrow v_m,$ $ m \ge 2,$
		either there is no arc between $v_0$ and $v_m$, or for all $1 \le i < j \le m$, the arc $v_i \rightarrow v_j$ exists. An induced subgraph on vertices $\{v_0, v_1, \ldots, v_k\}$ of an oriented graph is called a \emph{shortcut} if it is acyclic, non-transitive, and contains both the directed path $v_0 \rightarrow v_1 \rightarrow \cdots \rightarrow v_k$ and the arc $v_0 \rightarrow v_k$, referred to as the \emph{shortcutting edge}. Equivalently, an orientation is semi-transitive if and only if it is acyclic and contains no shortcuts.
	\end{definition}
	
	\begin{theorem}[\cite{halldorsson2011alternation}]\label{thm:wr_equiv_st}
		A graph is word-representable if and only if it admits a semi-transitive orientation.
	\end{theorem}
	This fundamental result provides a bridge between combinatorial properties of words and structural properties of graphs, allowing many problems on word-representability to be studied through orientations.
	\begin{remark}[\cite{kitaev2017comprehensive}]\label{rem:hereditary}
		The class of word-representable graphs is hereditary; that is, every induced subgraph of a word-representable graph is also word-representable.
	\end{remark}
	
	\medskip
	
	A word is said to be \emph{$k$-uniform} if each letter occurs exactly $k$ times. Correspondingly, a graph $G$ is \emph{$k$-word-representable} if it can be represented by a $k$-uniform word. The smallest integer $k$ for which a graph $G$ is $k$-word-representable is called the \emph{representation number} of $G$, denoted by $\mathcal{R}(G)$. We define
	\[
	\mathcal{R}_k := \{ G \mid \mathcal{R}(G) = k \}.
	\]
	
	\begin{theorem}[\cite{kitaev2008representable}]\label{thm:k_wordrep}
		A graph is word-representable if and only if it is $k$-word-representable for some $k$.
	\end{theorem}
	
	The class $\mathcal{R}_1$ consists precisely of complete graphs, whereas $\mathcal{R}_2$ coincides with the class of circle graphs excluding complete graphs. Recently, Fleischmann et~al.~\cite{fleischmann2024word} established that for $k \ge 3$, a word-representable graph is $k$-word-representable if and only if it admits a $k$-circle representation.
	
	\medskip
	
	The following useful results capture symmetry properties of word-representants.
	
	\begin{prop}[\cite{kitaev2015words}]\label{prop:rotation}
		Let $w = uv$ be a $k$-uniform word representing a graph $G$, where $u$ and $v$ are (possibly empty) subwords. Then the word $w' = vu$ also represents $G$.
	\end{prop}
	
	\begin{definition}
		The \emph{reverse} of a word $w = w_1w_2\ldots w_n$ is the word $r(w) = w_n\ldots w_2w_1$.
	\end{definition}
	
	\begin{prop}[\cite{kitaev2015words}]\label{prop:reverse}
		If $w$ is a word-representant of a graph $G$, then its reverse $r(w)$ also represents $G$.
	\end{prop}
	
		\begin{definition} A graph $G$ is \emph{permutationally representable} if it can be represented by a word of the form $p_1 \cdots p_k$, where each $p_i$ is a permutation. In this case, we say that $G$ is \emph{permutationally $k$-representable}, and the word $w$ is called a \textit{$k$-permutation word-representant} of $G$.
	\end{definition}
	
	\begin{theorem}[\cite{kitaev2008representable}]\label{neighbourhood} 
		If a graph $G$ is word-representable, then the neighborhood of each vertex in $G$ is permutationally representable.
	 \end{theorem}
	
	For a subset of letters $\{x_1, \ldots, x_m\}$, the \emph{restriction} of $w$ to these letters, denoted by $w|_{\{x_1, \ldots, x_m\}}$, is obtained by deleting from $w$ all letters not in $\{x_1, \ldots, x_m\}$. For example, if $w = 35423214$, then $w|_{\{1,2\}} = 221$. For further reading on word-representable graphs, we refer the reader to \cite{kitaev2015words, kitaev2017comprehensive}.	
	\section{Forbidden induced subgraph characterization of co-bipartite circle graphs}\label{s3}
	
	In this section, we establish the forbidden induced subgraph characterization of circle graphs within the class of co-bipartite graphs. Our main result demonstrates that the class of co-bipartite circle graphs coincides with the class of co-bipartite permutation graphs. Consequently, by applying Gallai’s forbidden induced subgraph characterization of comparability graphs~\cite{MR221974}, we obtain the corresponding characterization for co-bipartite circle graphs. The family of minimal forbidden induced subgraphs of permutation graphs within co-bipartite graphs, denoted by $\gc$, is illustrated in \textit{Figure}~\ref{cop}. Before presenting the proof of the main theorem, we recall several auxiliary lemmas and preliminaries on circle graphs that will be used in the argument.
	
	\begin{figure}[h]
	
		\begin{minipage}{0.45\textwidth}
			\centering
			\begin{tikzpicture}[line cap=round,line join=round,>=triangle 45,x=1cm,y=1cm, scale=2]
				\begin{scriptsize}
					\node[draw, circle, fill = qqzzcc] (k-1) at (1.2,3.25) {};
					\draw[color=qqzzcc] (1,3.4) node {2k-2};
					\node[draw, circle, fill = qqccqq] (3) at (1.8,3.25) {};
					\draw[color=qqccqq] (1.88,3.4) node {3};
					\node[draw, circle, fill = qqzzcc] (2) at (2.1,2.75) {};
					\draw[color=qqzzcc] (2.25,2.72) node {2};
					\node[draw, circle, fill = qqccqq] (1) at (1.8,2.25) {};
					\draw[color=qqccqq] (1.88,2.1) node {1};
					\node[draw, circle, fill = qqccqq] (k) at (0.9,2.75) {};
					\draw[color=qqccqq] (0.65,2.73) node {2k - 1};
					\draw[color=black] (1.5,3.4) node {...};
					\node[draw, circle, fill = qqzzcc] (k+1) at (1.2,2.25) {};
					\draw[color=qqzzcc] (1.1,2.1) node {2k};
					
					\draw (k-1) -- (3);
					\draw (1) -- (3);
					\draw (k) -- (3);
					\draw (k+1) -- (3);
					\draw (k-1) -- (1);
					\draw (k) -- (1);
					\draw (k-1) -- (2);
					\draw (k+1) -- (2);
					\draw (k) -- (2);
					\draw (k-1) -- (k+1);
				\end{scriptsize}
			\end{tikzpicture}\\{$\overline{C_{2k}}$ for $k \ge 3$}
		\end{minipage}
		\begin{minipage}{0.45\textwidth}
			\centering
			\begin{tikzpicture}[line cap=round,line join=round,>=triangle 45,x=1cm,y=1cm, scale=2]
				\begin{scriptsize}
					\node[draw, circle, fill = qqzzcc] (4) at (1.05,2.15) {};
					\draw[color=qqzzcc] (0.85,2.12) node {7};
					\node[draw, circle, fill = qqzzcc] (2) at (1.95,2.75) {};
					\draw[color=qqzzcc] (2.1,2.72) node {5};
					\node[draw, circle, fill = qqzzcc] (1) at (1.95,2.15) {};
					\draw[color=qqzzcc] (2.1,2.12) node {6};
					\node[draw, circle, fill = qqzzcc] (5) at (1.05,2.75) {};
					\draw[color=qqzzcc] (0.85,2.73) node {4};
					\node[draw, circle, fill = qqccqq] (c) at (1.95,3.35) {};
					\draw[color=qqccqq] (2.1,3.33) node {2};
					\node[draw, circle, fill = qqccqq] (d) at (1.5,3) {};
					\draw[color=qqccqq] (1.5,3.15) node {3};
					\node[draw, circle, fill = qqccqq] (b) at (1.05,3.35) {};
					\draw[color=qqccqq] (0.85,3.33) node {1};
					
					\draw (4) -- (1);
					\draw (2) -- (1);
					\draw (5) -- (1);
					\draw (4) -- (2);
					\draw (5) -- (2);
					\draw (4) -- (5);
					\draw (5) -- (d);
					\draw (2) -- (d);
					\draw (2) -- (c);
					\draw (5) -- (b);
					\draw (1) -- (d);
					\draw (c) -- (d);
					\draw (b) -- (d);
					\draw (b) -- (c);
				\end{scriptsize}
			\end{tikzpicture}\\{$G_1$}
		\end{minipage}
		
		\vspace{0.5cm}
		\noindent
		
		\begin{minipage}{0.44\textwidth}
			\centering
			\begin{tikzpicture}[line cap=round,line join=round,>=triangle 45,x=1cm,y=1cm, scale=2]
				\begin{scriptsize}
					\node[draw, circle, fill = qqzzcc] (4) at (1.5,2.3) {};
					\draw[color=qqzzcc] (1.35,2.55) node {4};
					\node[draw, circle, fill = qqzzcc] (2) at (1.8,2.15) {};
					\draw[color=qqccqq] (2.2,2.75) node {3};
					\node[draw, circle, fill = qqzzcc] (1) at (1.5,2.6) {};
					\draw[color=qqccqq] (1.5,3.35) node {1};
					\node[draw, circle, fill = qqzzcc] (5) at (1.2,2.15) {};
					\draw[color=qqccqq] (0.8,2.73) node {2};
					\node[draw, circle, fill = qqccqq] (c) at (2.05,2.75) {};
					\draw[color=qqzzcc] (1.05,2.15) node {5};
					\node[draw, circle, fill = qqccqq] (d) at (1.5,3.2) {};
					\draw[color=qqzzcc] (1.38,2.3) node {7};
					\node[draw, circle, fill = qqccqq] (b) at (0.95,2.75) {};
					\draw[color=qqzzcc] (1.95,2.15) node {6};
					
					\draw (4) -- (1);
					\draw (2) -- (1);
					\draw (5) -- (1);
					\draw (4) -- (2);
					\draw (5) -- (2);
					\draw (4) -- (5);
					\draw (c) -- (d);
					\draw (b) -- (d);
					\draw (b) -- (c);
					\draw (1) -- (c);
					\draw (2) -- (c);
					\draw (5) -- (b);
					\draw (1) -- (b);
				\end{scriptsize}
			\end{tikzpicture}\\{$G_2$}
		\end{minipage}
		\begin{minipage}{0.45\textwidth}
			\centering
			\begin{tikzpicture}[line cap=round,line join=round,>=triangle 45,x=1cm,y=1cm, scale=2]
				\begin{scriptsize}
						\node[draw, circle, fill = qqzzcc] (4) at (1.5,2.3) {};
					\draw[color=qqzzcc] (1.35,2.55) node {4};
					\node[draw, circle, fill = qqzzcc] (2) at (1.8,2.15) {};
					\draw[color=qqccqq] (2.2,2.75) node {3};
					\node[draw, circle, fill = qqzzcc] (1) at (1.5,2.6) {};
					\draw[color=qqccqq] (1.5,3.35) node {1};
					\node[draw, circle, fill = qqzzcc] (5) at (1.2,2.15) {};
					\draw[color=qqccqq] (0.8,2.73) node {2};
					\node[draw, circle, fill = qqccqq] (c) at (2.05,2.75) {};
					\draw[color=qqzzcc] (1.05,2.15) node {5};
					\node[draw, circle, fill = qqccqq] (d) at (1.5,3.2) {};
					\draw[color=qqzzcc] (1.38,2.3) node {7};
					\node[draw, circle, fill = qqccqq] (b) at (0.95,2.75) {};
					\draw[color=qqzzcc] (1.95,2.15) node {6};
					
					\draw (4) -- (1);
					\draw (2) -- (1);
					\draw (5) -- (1);
					\draw (4) -- (2);
					\draw (5) -- (2);
					\draw (4) -- (5);
					\draw (5) -- (d);
					\draw (2) -- (d);
					\draw (1) -- (c);
					\draw (2) -- (c);
					\draw (5) -- (b);
					\draw (1) -- (b);
					\draw (c) -- (d);
					\draw (b) -- (d);
					\draw (b) -- (c);
				\end{scriptsize}
			\end{tikzpicture}\\{$G_3$}
		\end{minipage}
		
		\vspace{0.25cm}\caption{Minimal forbidden induced subgraphs for co-bipartite permutation graphs, $\gc$, where vertices of same color induce a clique.}
		\label{cop}
	\end{figure}
	
	The following characterization of circle graphs using local complementation will be used in the proof of the results of this section.
	
	\begin{theorem}[\cite{MR1256586}]\label{circ}
		A graph $G$ is a circle graph if and only if $G$ has no $\ell$-reduction isomorphic to $W_5$, $W_7$, or $Y_6$.
	\end{theorem}
		\begin{figure}[h]
		\begin{minipage}{0.24\textwidth}
			\centering
			\begin{tikzpicture}[line cap=round,line join=round,>=triangle 45,x=1cm,y=1cm, scale=2]
				\begin{scriptsize}
					\node[draw, circle, fill = qqzzcc] (2) at (1.8,2.15) {};
					\node[draw, circle, fill = qqzzcc] (1) at (1.5,2.6) {};
					\node[draw, circle, fill = qqzzcc] (5) at (1.2,2.15) {};
					\node[draw, circle, fill = qqzzcc] (c) at (2.05,2.75) {};
					\node[draw, circle, fill = qqzzcc] (d) at (1.5,3.2) {};
					\node[draw, circle, fill = qqzzcc] (b) at (0.95,2.75) {};
					
					\draw (2) -- (1);
					\draw (5) -- (1);
					\draw (5) -- (2);
					\draw (c) -- (d);
					\draw (b) -- (d);
					\draw (1) -- (c);
					\draw (2) -- (c);
					\draw (5) -- (b);
					\draw (1) -- (b);
					\draw (1) -- (d);
				\end{scriptsize}
			\end{tikzpicture}\\{$W_5$}
		\end{minipage}
		\hfill
		\begin{minipage}{0.24\textwidth}
			\centering
			\begin{tikzpicture}[line cap=round,line join=round,>=triangle 45,x=1cm,y=1cm, scale=2]
				\begin{scriptsize}
					\node[draw, circle, fill = qqzzcc] (2) at (1.8,2.15) {};
					\node[draw, circle, fill = qqzzcc] (3) at (1.95,2.45) {};
					\node[draw, circle, fill = qqzzcc] (4) at (1.05,2.45) {};
					\node[draw, circle, fill = qqzzcc] (1) at (1.5,2.6) {};
					\node[draw, circle, fill = qqzzcc] (5) at (1.2,2.15) {};
					\node[draw, circle, fill = qqzzcc] (c) at (2.05,2.75) {};
					\node[draw, circle, fill = qqzzcc] (d) at (1.5,3.2) {};
					\node[draw, circle, fill = qqzzcc] (b) at (0.95,2.75) {};
					
					\draw (4) -- (1);
					\draw (2) -- (1);
					\draw (3) -- (1);
					\draw (5) -- (1);
					\draw (5) -- (2);
						\draw (3) -- (2);
							\draw (3) -- (c);
					\draw (1) -- (c);
					\draw (4) -- (b);
					\draw (5) -- (4);
					\draw (1) -- (b);
					\draw (1) -- (d);
					\draw (c) -- (d);
					\draw (b) -- (d);
				\end{scriptsize}
			\end{tikzpicture}\\{$W_7$}
		\end{minipage}
		\hfill
		\begin{minipage}{0.24\textwidth}
			\centering
			\begin{tikzpicture}[line cap=round,line join=round,>=triangle 45,x=1cm,y=1cm, scale=2]
				\begin{scriptsize}
					\node[draw, circle, fill = qqzzcc] (4) at (1.5,2) {};
					\node[draw, circle, fill = qqzzcc] (2) at (2.05,2.2) {};
					\node[draw, circle, fill = qqzzcc] (1) at (1.5,2.5) {};
					\node[draw, circle, fill = qqzzcc] (5) at (0.95,2.2) {};
					\node[draw, circle, fill = qqzzcc] (c) at (2.05,2.75) {};
					\node[draw, circle, fill = qqzzcc] (d) at (1.5,3) {};
					\node[draw, circle, fill = qqzzcc] (b) at (0.95,2.75) {};
					
					\draw (4) -- (1);
					\draw (4) -- (2);
					\draw (4) -- (5);
					\draw (1) -- (c);
					\draw (2) -- (c);
					\draw (5) -- (b);
					\draw (1) -- (b);
					\draw (c) -- (d);
					\draw (b) -- (d);
				\end{scriptsize}
			\end{tikzpicture}\\{$Y_6$}
		\end{minipage}
		\vspace{0.25cm}\caption{Circle graph obstructions.}
		\label{circle}
	\end{figure}
	
			The following result shows that no graph in $\gc$ is a circle graph.
			The graphs $G_1$, $G_2$, and $G_3$ are excluded by noting that each admits an $l$-reduction isomorphic to one of the graphs $W_5$ or $Y_6$.
			For $G \cong \overline{C_{2k}}$, we argue by contradiction: assuming the existence of a $2$-uniform word-representant, we restrict it to a suitable closed neighborhood and obtain a permutation graph.
			We then show that the corresponding permutation word is essentially unique.
			This rigidity enforces strict alternation constraints, leaving no admissible positions for the remaining vertices.
			The resulting contradiction implies that $\overline{C_{2k}}$, and hence every graph in $\gc$, is not a circle graph.

	\begin{lemma}\label{gcp}
		If $G \in \gc$, then $G$ is not a circle graph.
	\end{lemma}
	\begin{proof}
		Let $G$ be a graph in $\gc$. Consider the vertex set of $G$ as provided in \textit{Figure~}\ref{cop}.\\
		\textit{Case~1 (local complementation of $G_1$ leading to $W_5$).}		
		 Suppose $G \cong G_1$. It follows that $G \ast 7~6$ contains an induced subgraph isomorphic to $W_5$.\\
		\textit{Case~2 (local complementation of $G_2$ leading to $Y_6$).}		
		Suppose $G \cong G_2$. It follows that $G \ast 7~1$ contains an induces subgraph isomorphic to $Y_6$.\\
		\textit{Case~3 (local complementation of $G_3$ leading to $Y_6$).}
		Suppose $G \cong G_3$. It follows that $G \ast 1~2~3$ contains an induces subgraph isomorphic to $Y_6$.
			\begin{figure}[h]
			
			\begin{minipage}{0.32\textwidth}
				\centering
				\begin{tikzpicture}[line cap=round,line join=round,>=triangle 45,x=1cm,y=1cm, scale=2]
					\begin{scriptsize}
						\node[draw, circle, fill = qqzzcc] (4) at (1.05,2.15) {};
						\draw[color=qqzzcc] (0.85,2.12) node {7};
						\node[draw, circle, fill = qqzzcc] (2) at (1.95,2.75) {};
						\draw[color=qqzzcc] (2.1,2.72) node {5};
						\node[draw, circle, fill = qqzzcc] (1) at (1.95,2.15) {};
						\draw[color=qqzzcc] (2.1,2.12) node {6};
						\node[draw, circle, fill = qqzzcc] (5) at (1.05,2.75) {};
						\draw[color=qqzzcc] (0.85,2.73) node {4};
						\node[draw, circle, fill = qqzzcc] (c) at (1.95,3.35) {};
						\draw[color=qqzzcc] (2.1,3.33) node {2};
						\node[draw, circle, fill = qqzzcc] (d) at (1.5,3) {};
						\draw[color=qqzzcc] (1.5,3.15) node {3};
						\node[draw, circle, fill = qqzzcc] (b) at (1.05,3.35) {};
						\draw[color=qqzzcc] (0.85,3.33) node {1};
						
						\draw (4) -- (1);
						\draw (d) -- (4);
						\draw (4) -- (2);
						\draw (4) -- (5);
						\draw (5) -- (d);
						\draw (2) -- (d);
						\draw (2) -- (c);
						\draw (5) -- (b);
						\draw (1) -- (d);
						\draw (c) -- (d);
						\draw (b) -- (d);
						\draw (b) -- (c);
					\end{scriptsize}
				\end{tikzpicture}\\{$G_1 \ast 7~6$}
			\end{minipage}
			\begin{minipage}{0.32\textwidth}
				\centering
				\begin{tikzpicture}[line cap=round,line join=round,>=triangle 45,x=1cm,y=1cm, scale=2]
					\begin{scriptsize}
						\node[draw, circle, fill = qqzzcc] (4) at (1.5,2.3) {};
						\draw[color=qqzzcc] (1.35,2.55) node {4};
						\node[draw, circle, fill = qqzzcc] (2) at (1.8,2.15) {};
						\draw[color=qqzzcc] (2.2,2.75) node {3};
						\node[draw, circle, fill = qqzzcc] (1) at (1.5,2.6) {};
						\draw[color=qqzzcc] (1.5,3.35) node {1};
						\node[draw, circle, fill = qqzzcc] (5) at (1.2,2.15) {};
						\draw[color=qqzzcc] (0.8,2.73) node {2};
						\node[draw, circle, fill = qqzzcc] (c) at (2.05,2.75) {};
						\draw[color=qqzzcc] (1.05,2.15) node {5};
						\node[draw, circle, fill = qqzzcc] (d) at (1.5,3.2) {};
						\draw[color=qqzzcc] (1.38,2.3) node {7};
						\node[draw, circle, fill = qqzzcc] (b) at (0.95,2.75) {};
						\draw[color=qqzzcc] (1.95,2.15) node {6};
						
						\draw (4) -- (1);
						\draw (4) -- (2);
						\draw (4) -- (5);
						\draw (c) -- (d);
						\draw (b) -- (d);
						\draw (1) -- (c);
						\draw (2) -- (c);
						\draw (5) -- (b);
						\draw (1) -- (b);
					\end{scriptsize}
				\end{tikzpicture}\\{$G_2 \ast 7~1$}
			\end{minipage}
			\begin{minipage}{0.32\textwidth}
				\centering
				\begin{tikzpicture}[line cap=round,line join=round,>=triangle 45,x=1cm,y=1cm, scale=2]
					\begin{scriptsize}
						\node[draw, circle, fill = qqzzcc] (4) at (1.5,2.3) {};
						\draw[color=qqzzcc] (1.35,2.55) node {4};
						\node[draw, circle, fill = qqzzcc] (2) at (1.8,2.15) {};
						\draw[color=qqzzcc] (2.2,2.75) node {3};
						\node[draw, circle, fill = qqzzcc] (1) at (1.5,2.6) {};
						\draw[color=qqzzcc] (1.5,3.35) node {1};
						\node[draw, circle, fill = qqzzcc] (5) at (1.2,2.15) {};
						\draw[color=qqzzcc] (0.8,2.73) node {2};
						\node[draw, circle, fill = qqzzcc] (c) at (2.05,2.75) {};
						\draw[color=qqzzcc] (1.05,2.15) node {6};
						\node[draw, circle, fill = qqzzcc] (d) at (1.5,3.2) {};
						\draw[color=qqzzcc] (1.38,2.3) node {7};
						\node[draw, circle, fill = qqzzcc] (b) at (0.95,2.75) {};
						\draw[color=qqzzcc] (1.95,2.15) node {5};
						
						\draw (4) -- (1);
						\draw (4) -- (2);
						\draw (4) -- (5);
						\draw (c) -- (d);
						\draw (b) -- (d);
						\draw (1) -- (c);
						\draw (2) -- (c);
						\draw (5) -- (b);
						\draw (1) -- (b);
					\end{scriptsize}
				\end{tikzpicture}\\{$G_3 \ast 1~2~3$}
			\end{minipage}
			\caption{Local complementations of $G_1$, $G_2$, and $G_3$ yield $W_5$ or $Y_6$.}
		\end{figure}
		
		Therefore, for the above three cases, by \textit{Theorem~\ref{circ}}, it follows that $G$ is not a circle graph. Hence, suppose that $G \cong \overline{C_{2k}}$ for some $k \ge 3$. Assume, for the sake of contradiction, that $G$ is a circle graph. This implies that there exists a $2$-uniform word-representant of $G$, say $w_G$. Consider the graph $H \cong G[N_G[1]]$. By \textit{Theorem~\ref{neighbourhood}} and the hereditary nature of circle graphs, it can be seen that $H$ is a permutation graph. Let $w_H$ denote the subword of $w_G$ obtained by restricting its alphabet to the vertex set of $H$. Then, it is evident that $w_H$ is a $2$-permutation word-representant of $H$. By \textit{Proposition~\ref{prop:rotation}}, we may assume that the first letter of $w_H$ is $3$. We claim that 
		\[
		w_H = 3~5~4~7~6~\cdots~(2k-1)~(2k-2)~1~4~3~6~5~\cdots~(2k-2)~(2k-3)~(2k-1)~1
		\]
		is the unique 2-permutation word-representant of $H$, up to cyclic shifts and reversals. 
		
		The proof of the claim proceeds by induction on $k$. For the base case, consider $k = 3$. Since the first letter of $w_H$ is $3$, the two occurrences of $4$ in $w_H$ must appear immediately after the first occurrence and immediately before the second occurrence of $3$, respectively. Similarly, the two occurrences of $5$ in $w_H$ must appear between the first occurrences and after the second occurrences of $3$ and $4$, respectively. Finally, the two occurrences of $1$ are placed after the first and second occurrences of $3$, $4$, and $5$, respectively. Therefore,
		\[
		w_H = 3~5~4~1~4~3~5~1.
		\]
		\noindent
		This establishes the base case of the induction. Assume that the claim holds for $k = i - 1$ for some $i > 4$. Then, we have
		\[
		w_H = 3~5~4~7~6~\cdots~(2i-3)~(2i-4)~1~4~3~6~5~\cdots~(2i-4)~(2i-5)~(2i-3)~1.
		\]
		Now, consider $k = i$. It suffices to show that there exists a unique way, up to cyclic shifts and reversals, to place the vertices $2i-2$ and $2i-1$ in $w_H$. The vertex $2i-2$ is adjacent to all vertices except $2i-3$ and $2i-1$. Therefore, if the first occurrence of $2i-2$ is placed before the first occurrence of $2i-3$, then, to maintain the non-alternation between $2i-2$ and $2i-3$, the second occurrence of $2i-2$ must be placed after the second occurrence of $2i-3$. However, this placement would result in a non-alternation between $2i-2$ and $2i-4$, contradicting the existence of a $2$-permutation word-representant of $H$. Hence, the unique placement of vertex $2i-2$ is as follows:
		\[
		w_H = 3~5~4~7~6~\cdots~(2i-3)~(2i-4)~(2i-2)~1~4~3~6~5~\cdots~(2i-4)~(2i-5)~(2i-2)~(2i-3)~1.
		\]
		\noindent
		Similarly, the vertex $2i-1$ is adjacent to all vertices except $2i-2$ and $2i$. Therefore, if the first occurrence of $2i-1$ is placed after the first occurrence of $2i-2$, then, to maintain the non-alternation between $2i-1$ and $2i-2$, the second occurrence of $2i-1$ must be placed before the second occurrence of $2i-2$. However, this placement would result in a non-alternation between $2i-1$ and $2i-3$, contradicting the existence of a $2$-permutation word-representant of $H$. Hence, the first occurrence of $2i-1$ must be placed before the first occurrence of $2i-2$. 
		
		Furthermore, if the first occurrence of $2i-1$ is placed before the first occurrence of $2i-4$, then, to maintain the alternation between $2i-1$ and $2i-4$, the second occurrence of $2i-1$ must be placed before the second occurrence of $2i-4$. This, however, would lead to an alternation between $2i-1$ and $2i-2$, again a contradiction. Therefore, the unique placement of the vertex $2i-1$ is as follows:
		\[
		w_H = 3~5~4~7~6~\cdots~(2i-4)~(2i-1)~(2i-2)~1~4~3~6~5~\cdots~(2i-2)~(2i-3)~(2i-1)~1.
		\]
		\noindent
		Hence, the placement of the vertices $2i-2$ and $2i-1$ in $w_H$ is unique up to cyclic shifts and reversals, thereby proving the claim of the induction hypothesis. 
		
		Now consider $w_G$. Since $w_H$ is a subword of $w_G$, we can obtain $w_G$ by appropriately placing the vertices $2$ and $2k$ in $w_H$. The vertex $2$ is adjacent to all vertices except $1$ and $3$. Suppose the first occurrence of $2$ is placed before the first occurrence of $3$. Then, to maintain the non-alternation between $2$ and $3$, the second occurrence of $2$ must appear after the second occurrence of $3$. Furthermore, to preserve the non-alternation between $2$ and $1$, the second occurrence of $2$ must also appear after the second occurrence of $1$. However, such a placement would create a non-alternation between $2$ and another vertex, thereby contradicting the assumption that $G$ admits a $2$-word-representant. Hence, the first occurrence of $2$ must be placed after the first occurrence of $3$. 
		
		If, instead, the first occurrence of $2$ is placed after the first occurrence of $5$, then to maintain the alternation between $2$ and $5$, the second occurrence of $2$ must also appear after the second occurrence of $5$. This, however, results in a non-alternation between $2$ and $4$, contradicting the existence of a $2$-word-representant of $G$. Therefore, the first occurrence of $2$ must lie between the first occurrences of $3$ and $5$. Consequently, to maintain the non-alternation between $2$ and $1$, and alternation between $2$ and all other vertices except $1$ and $3$, the second occurrence of $2$ must be placed between the first occurrence of $2i-2$ and $1$. Hence,
		\[
		w_G|_{V(G)\setminus\{2k\}} = 3~2~5~4~7~6~\cdots~(2k-1)~(2k-2)~2~1~4~3~6~5~\cdots~(2k-3)~(2k-1)~1.
		\]
		
		The vertex $2k$ is adjacent to all vertices except $2k-1$ and $1$. Suppose that the first occurrence of $2k$ appears before the first occurrence of $2k-1$. Then, to maintain the non-alternation between $2k$ and $2k-1$, either the second occurrence of $2k$ must appear before the first occurrence of $2k-1$, or after the second occurrence of $2k-1$. However, in both cases, a non-alternation arises between $2k$ and $2k-2$, contradicting the existence of a $2$-word-representant of $G$. Therefore, assume that the first occurrence of $2k$ appears after the first occurrence of $2k-1$. In this case, to preserve the non-alternation between $2k$ and $2k-1$, and the alternation between $2k$ and $2k-3$, the second occurrence of $2k$ must occur between the second occurrences of $2k-3$ and $2k-1$. Furthermore, to maintain the non-alternation between $1$ and $2k$, the first occurrence of $2k$ must appear after the first occurrence of $1$. This, however, introduces a non-alternation between $2k$ and $2$, again leading to a contradiction. Hence, vertex $2k$ cannot appear anywhere in $w_G$, contradicting the assumption that $G$ is $2$-word-representable. Therefore, $G$ is not a circle graph. This completes the proof of the lemma.
	\end{proof}
	\begin{theorem}
		Let $G$ be a co-bipartite graph. Then $G$ is a circle graph if and only if it is a permutation graph.
	\end{theorem}
	\begin{proof}
		Since permutation graphs form a subclass of circle graphs, the converse follows trivially. Thus, it suffices to prove the forward direction. Suppose that $G$ is a circle graph. For the sake of contradiction, assume that $G$ is not a permutation graph. Then, by the characterization of permutation graphs within co-bipartite graphs in terms of minimal forbidden induced subgraphs, $G$ contains some graph in $\gc$. However, by \textit{Lemma}~\ref{gcp} and the hereditary property of circle graphs, it follows that $G$ is not a circle graph, a contradiction.
	\end{proof}
	\begin{corollary}\label{fisccirc}
		Let $G$ be a co-bipartite graph. Then $G$ is a circle graph if and only if it contains no graph in $\gc$ as an induced subgraph.
	\end{corollary}
	Having established a forbidden induced subgraph characterization of co-bipartite circle graphs, we now turn our attention to the broader class of word-representable co-bipartite graphs. Instead of working directly with graph orientations, we adopt a matrix-theoretic perspective, which enables us to leverage existing characterizations and algorithms for circular-ones-type properties.
	\section{Forbidden induced subgraph characterization of word-representable co-bipartite graphs}\label{s4}
	
	For a co-bipartite graph $G = (X, Y)$, we define $M(G) = \{m_{ij}\}$ to be a $(0,1)$-matrix whose rows correspond to the partition $X$ and whose columns correspond to the partition $Y$, where $m_{ij} = 1$ if and only if the corresponding row vertex and column vertex are adjacent in $G$.	
	In this section, we establish a connection between the circularly compatible ones property of matrices, introduced by Tucker in~\cite{MR2618994}, and semi-transitive co-bipartite graphs, thereby providing a forbidden induced subgraph characterization of the latter. Since the characterization of word-representable co-bipartite graphs given by Das and Ramesh in~\cite{das2025representation} involves monotonicity, we are naturally motivated to consider the circularly compatible ones property, which also incorporates a monotonicity condition.
	
	In~\cite{mdsafe2}, Safe provided a minimal forbidden submatrix characterization and a linear-time recognition algorithm for the circularly compatible ones property. 
	Consequently, by establishing a connection between semi-transitive co-bipartite graphs and the circularly compatible ones property of matrices, we obtain a minimal forbidden induced subgraph characterization for semi-transitive co-bipartite graphs. Before presenting the main result, we summarize the preliminaries concerning the circularly compatible ones property. Throughout this section, for any two vertices $x$ and $y$ in an oriented graph, we write $x \rightarrow y$ if there exists an arc between $x$ and $y$ in the given orientation.
	
	\begin{definition}
		Let $\preccurlyeq_x$ be a linear order on a set $X$. For elements $a,b \in X$ with $a \preccurlyeq_x b$, the \emph{linear interval of $\preccurlyeq_x$ with endpoints $a$ and $b$}, denoted by $[a,b]_{\preccurlyeq_x}$, is defined as the set $\{x \in X : a \preccurlyeq_x x \preccurlyeq_x b\}$. A \emph{linear interval of $\preccurlyeq_x$} is either the empty set or an interval of the form $[a,b]_{\leq_X}$ for some $a,b \in X$ satisfying $a \preccurlyeq_x b$. A sequence $a_1a_2\ldots a_k$ is said to be \emph{monotone} on $X$ if $a_1,a_2,\ldots,a_k \in X$ and $a_1 \preccurlyeq_x a_2 \preccurlyeq_x \cdots \preccurlyeq_x a_k$.
	\end{definition}
	\begin{definition}
		A matrix $M$ has the consecutive-ones property for rows if there is a linear order $\preccurlyeq_C$ of the columns of $M$ such that each row of $M$ is a linear interval of $\preccurlyeq_C$, and $\preccurlyeq_C$ is called a consecutive-ones order. Analogous definition apply to the columns of $M$. If no mention is made for rows or columns, we mean the corresponding property for the rows.
	\end{definition}
	\begin{definition}
		Let $\preccurlyeq_x$ be a linear order on a set $X$, and let $a,b \in X$. 
		The \emph{circular interval of $\preccurlyeq_x$ with endpoints $a$ and $b$}, denoted by $[a,b]_{\preccurlyeq_x}$, is defined as follows:  
		if $a \preccurlyeq_x b$, then 
		\[
		[a,b]_{\preccurlyeq_x} = \{x \in X \mid a \preccurlyeq_x x \preccurlyeq_x b\},
		\] 
		while if $b \prec_X a$, then 
		\[
		[a,b]_{\preccurlyeq_x} = \{x \in X \mid x \preccurlyeq_x b \text{ or } a \preccurlyeq_x x\}.
		\] 
		A \emph{circular interval of $\preccurlyeq_x$} is either the empty set or a set of the form $[a,b]_{\preccurlyeq_x}$ for some $a,b \in X$. A sequence $a_1a_2\ldots a_k$ is said to be \emph{circularly monotone with respect to $\preccurlyeq_x$} if $a_1,a_2,\ldots,a_k \in X$ and the relation $a_i \preccurlyeq_x a_{i+1}$ holds for all but at most one index $i \in \{1,2,\ldots,k\}$, where $a_{k+1}$ is identified with $a_1$.
	\end{definition}
	\begin{definition}\label{cir}
		A matrix $M$ has the circular-ones property for rows if there is a linear order $\preccurlyeq_C$ of the columns of $M$ such that each row of $M$ is a circular interval of $\preccurlyeq_C$, and $\preccurlyeq_C$ is called a circular-ones order. Analogous definition apply to the columns of $M$. If no mention is made for rows or columns, we mean the corresponding property for the rows.
	\end{definition}
	The above definitions introduce linear and circular intervals with respect to a given linear order, together with the corresponding notions of monotonicity. 
	These concepts form the basis for defining the consecutive-ones and circular-ones properties of a matrix.	
	
	\begin{theorem}[\cite{MR309810}]\label{c1cp}
		Let $M$ be a matrix, and let $M'$ denote any matrix obtained from $M$ by complementing some of its rows so that $M'$ contains at least one column consisting entirely of zeros. Then, $M$ has the circular-ones property if and only if $M'$ has the consecutive-ones property.
	\end{theorem}
	\begin{theorem}[\cite{mdsafe1}]\label{circR} A matrix $M$ has the circular-ones property if and only if $M$ contains no matrix in the set $\ForbRow$ as a configuration. The corresponding set of forbidden submatrices is
		\[ \ForbRow=\{a\miop\MIast k\colon\,k\geq 3\mbox{ and }a\in A_k\}\cup\{\MIV,\overline{\MIV},\MVast,\overline{\MVast}\}, \]
		where $\MIast k$ and $\MVast$ denote $(\MI k)^*$ and $(\MV)^*$, $A_3=\{000,111\}$ and, for each $k\geq 4$, $A_k$ is the set of all binary bracelets of length $k$. Notice that $001$ and $011$ are binary bracelets of length $3$ but do not belong to $A_3$. 
	\end{theorem}
	The following definitions present two strengthened variants of the circular-ones property for matrices, originally introduced by Safe in \cite{mdsafe2}.	
	\begin{definition}
		A matrix $M$ is said to have the \emph{$D$-circular property} if there exists a linear order $\preccurlyeq_c$ of its columns such that every row of $M$ is a circular interval with respect to $\preccurlyeq_c$, and for any two rows $r$ and $s$ of $M$, the set difference $s - r$ also forms a circular interval of $\preccurlyeq_c$.  
		In this case, $\preccurlyeq_c$ is called a \emph{$D$-circular order} of $M$.
	\end{definition}
	\begin{definition}
		A matrix $M$ is said to have the \emph{circularly compatible ones property} if there exists a biorder $(\preccurlyeq_r,\preccurlyeq_c)$ such that:  
		\begin{enumerate}[(i)]
			\item each row of $M$ forms a circular interval with respect to $\preccurlyeq_c$;  
			\item each column of $M$ forms a circular interval with respect to $\preccurlyeq_r$; and  
			\item if $r_1,r_2,\ldots,r_p$ denote all nontrivial rows of $M$ listed in ascending order of $\preccurlyeq_r$, and if $r_i = [d_i,e_i]_{\preccurlyeq_c}$ for each $i \in \{1,2,\ldots,p\}$, then both sequences $d_1d_2\ldots d_p$ and $e_1e_2\ldots e_p$ are circularly monotone with respect to $\preccurlyeq_c$.  
		\end{enumerate}  
		In this case, the pair $(\preccurlyeq_r,\preccurlyeq_c)$ is called a \emph{circularly compatible ones biorder}.
		
	\end{definition}
	
	In~\cite{MR3065109}, Basu et al.\ proved that, for a matrix $M$ without trivial rows, the $D$-circular property coincides with the monotone circular property. 
	We now state the definitions and results related to the monotone circular property.
	
	\begin{definition}\label{ui}
		Let $X = \{x_1, x_2, \ldots, x_k\}$ with $x_1 \preccurlyeq_x x_2 \preccurlyeq_x \cdots \preccurlyeq_x x_k$.  
		Define the set $X^+$ to be the set $\{x_1, x_2, \ldots, x_k, x_1^+, x_2^+, \ldots, x_k^+\}$, and let $\preccurlyeq_x^+$ be the linear order on $X^+$ given by  
		\[
		x_1 \preccurlyeq_x^+ x_2 \preccurlyeq_x^+ \cdots \preccurlyeq_x^+ x_k 
		\preccurlyeq_x^+ x_1^+ \preccurlyeq_x^+ x_2^+ \preccurlyeq_x^+ \cdots \preccurlyeq_x^+ x_k^+.
		\]  
		For each nontrivial circular interval $[a,b]_{\preccurlyeq_x}$, we define its \emph{unwrapped interval relative to $\preccurlyeq_x$}, denoted by $[a,b]_{\preccurlyeq_x}^+$, as the linear interval $[a,c]_{\preccurlyeq_x^+}$, where $c = b$ if $a \preccurlyeq_x b$, and $c = b^+$ if $b \prec_x a$.
	\end{definition}
	\begin{definition}\label{mc}
		Let $M$ be a matrix and let $(\preccurlyeq_r, \preccurlyeq_c)$ be a biorder of $M$.  
		Let $r_1, r_2, \ldots, r_p$ denote all nontrivial rows of $M$ listed in ascending order with respect to $\preccurlyeq_r$.  
		If $\preccurlyeq_c$ is a circular-ones order of $M$, $r_i = [d_i, e_i]_{\preccurlyeq_c}$ for each $i \in [p]$, and $[d_i, f_i]_{\preccurlyeq_c}^+$ denotes the unwrapped interval of $[d_i, e_i]_{\preccurlyeq_c}$ with respect to $\preccurlyeq_c$, we say that $M$ has \emph{monotone left endpoints} with respect to $(\preccurlyeq_r, \preccurlyeq_c)$ if the sequence $d_1 d_2 \ldots d_p$ is monotone with respect to $\preccurlyeq_c$, and we say that $M$ has \emph{monotone unwrapped right endpoints} with respect to $(\preccurlyeq_r, \preccurlyeq_c)$ if the sequence $f_1 f_2 \ldots f_p$ is monotone with respect to $\preccurlyeq_c^+$.  
		
		If $M$ has no trivial rows, a \emph{monotone circular biorder} of $M$ is a biorder $(\preccurlyeq_r, \preccurlyeq_c)$ satisfying all of the following conditions:
		\begin{enumerate}[(i)]
			\item $M$ has monotone left endpoints $d_1 d_2 \ldots d_p$ with respect to $(\preccurlyeq_r, \preccurlyeq_c)$;  
			\item $M$ has circularly monotone unwrapped right endpoints $f_1 f_2 \ldots f_p$ with respect to $(\preccurlyeq_r, \preccurlyeq_c)$;  
			\item either $f_1 = e_1^+$, or both $f_1 = e_1$ and $f_p \preccurlyeq_c^+ e_1^+$.  
		\end{enumerate}
		A matrix $M$ having no trivial rows has the \emph{monotone circular property} if it admits a monotone circular biorder. 
	\end{definition}
	
	The following result from~\cite{MR3065109}, which will be used in the proof of the main theorem of this section, relates the $D$-circular property to the monotone circular property.
	
	\begin{theorem}[\cite{MR3065109}]\label{mco}
		If $M$ is a matrix with no trivial rows, the following statements are equivalent:
		\begin{enumerate}[(i)]
			\item $M$ has the $D$-circular property;  
			\item $M$ has the monotone circular property.
		\end{enumerate}
	\end{theorem}
	
	\begin{definition}
		A matrix $M$ is said to have the \emph{doubly $D$-circular property} if both $M$ and its transpose $M^t$ possess the $D$-circular property.
	\end{definition}
	
	\begin{definition}
		We denote
		\begin{align*}
			\FCCO^\infty=\FCCO\cup\bigcup_{k=3}^\infty\{\MIast k,\overline{\MIast k},\MIast k\trans,\overline{\MIast k}\trans\},
		\end{align*}
		where
		\[   \FCCO=\{\Ztwo,\Ztthree,\Zfour,\Zfive,\CoZtwo,\CoZfour,
		(\Ztwo)\trans,(\Ztthree)\trans,(\Zfour)\trans,\ZfiveTrans,(\CoZtwo)\trans,(\CoZfour)\trans\}. \] (See Figure~\ref{f1})
	\end{definition}
	
	\begin{figure}[h]
		\ffigbox[\textwidth]{%
			
			\begin{subfloatrow}
				
				\subfloat[$\Ztwo$]{
					\begin{math}
						\left(\begin{array}{cccc}
							1 & 0 & 0 & 0\\
							1 & 1 & 0 & 0\\
							1 & 1 & 1 & 0\\
							0 & 1 & 0 & 0
						\end{array}\right)
					\end{math}
				}\quad
				\subfloat[$\Ztthree$]{
					\begin{math}
						\left(\begin{array}{cccc}
							1 & 0 & 0 & 0\\
							1 & 1 & 0 & 0\\
							1 & 1 & 1 & 0\\
							1 & 0 & 1 & 0
						\end{array}\right)
					\end{math}
				}\quad
				\subfloat[$\Zfour$]{
					\begin{math}
						\left(\begin{array}{ccccc}
							1 & 1 & 1 & 0 & 0\\
							0 & 1 & 1 & 1 & 0\\
							0 & 0 & 1 & 0 & 0
						\end{array}\right)
					\end{math}
				}
			\end{subfloatrow}
			\begin{subfloatrow}
				\subfloat[$\Zfive$]{
					\begin{math}
						\left(\begin{array}{cccc}
							1 & 0 & 0 & 1\\
							1 & 1 & 0 & 0\\
							1 & 1 & 1 & 0\\
							0 & 1 & 0 & 0
						\end{array}\right)
					\end{math}
				}\quad
				\subfloat[$\MIast k$ for each $k\geq 3$, where omitted entries are $0$'s]{
					\begin{math}
						\left(\begin{array}{cccccc}
							1 & 1 &        &       &   & 0\\
							& 1 & 1      &       &   & 0\\
							&   & \ddots & \ddots &  & \vdots \\
							&   &        &     1 & 1 & 0\\
							1 & 0 & \cdots &     0 & 1 & 0
						\end{array}\right)
					\end{math}
				}
			\end{subfloatrow}
		}{\caption{Some matrices in $\FCCO^\infty$}
			\label{f1}}
	\end{figure}
	The following result from~\cite{mdsafe2}, which will be used in the proof of the main theorem of this section, provides a minimal forbidden submatrix characterization of matrices possessing the circularly compatible ones property.
	
	\begin{theorem}[\cite{mdsafe2}]\label{cco}
		For any matrix $M$, the following statements are equivalent:
		\begin{enumerate}[(i)]
			\item $M$ has the circularly compatible ones property;  
			\item $M$ contains no matrix from $\FCCO^\infty$ as a configuration;  
			\item $M$ has the circular-ones property for both rows and columns, and contains no matrix from $\FCCO$ as a configuration;  
			\item $M$ has the doubly $D$-circular property.  
		\end{enumerate}
	\end{theorem}
	\begin{theorem}[\cite{mdsafe2}]\label{lta} 
		There exists a linear-time algorithm that, for any given matrix $M$, either produces a circularly compatible ones biorder of $M$ or identifies a matrix in $\FCCO^\infty$ that occurs in $M$ as a configuration.
	\end{theorem}
	For a $(0,1)$-matrix $M$ of size $m \times n$, we define $CG(M) = (X, Y)$ to be the co-bipartite graph in which the bipartition $(X, Y)$ corresponds to the rows and columns of $M$, respectively. We take $X = \{1, \dots, m\}$ and $Y = \{1, \dots, n\}$. A vertex $i \in X$ is adjacent to a vertex $j \in Y$ if and only if the corresponding entry of $M$ is $1$. By definition, $CG(M)$ is isomorphic to $CG(M\trans)$.
	\begin{lemma}\label{cgf}
		If $F \in \FCCO^\infty$, then $CG(F)$ is not semi-transitive.
	\end{lemma}
	\begin{proof}
	Let $F \in \FCCO^\infty$. Since $CG(F)$ is isomorphic to $CG(F\trans)$, it suffices to prove that $CG(F)$ is not semi-transitive for $F \in \{\MIast{k}, \overline{\MIast{k}} \colon k \ge 3\} \cup \{\Ztwo, \Ztthree, \Zfour, \Zfive, \CoZtwo, \CoZfour\}$. If $F = \overline{\MIast{k}}$ for some $k \ge 3$, then one can verify that $CG(F) \cong \overline{C_{2k}} \times K_1$. When $F \in \{\CoZtwo, \CoZfour\}$, we have $CG(F) \cong G_1 \times K_1$, and when $F \in \{\Ztthree, \Zfour\}$, we have $CG(F) \cong G_2 \times K_1$. In all these cases, $CG(F) \cong H \times K_1$ for some $H \in \gc$. Hence, by \textit{Theorem~\ref{neighbourhood}}, it follows that $CG(F)$ is not semi-transitive. Finally, if $F \in \{\Ztwo, \Zfive\}$, then by \textit{Theorem~9} of~\cite{chen2025word}, $CG(F)$ is not semi-transitive.
	
	Suppose that $F = \MIast{k}$ for some $k \ge 3$. Assume, for the sake of contradiction, that $CG(F)$ admits a semi-transitive orientation. Let the rows and columns of $F$ be denoted by $\{r_1, r_2, \ldots, r_k\}$ and $\{c_1, c_2, \ldots, c_{k+1}\}$, respectively. Define linear orders $\preccurlyeq_r$ and $\preccurlyeq_c$ on the rows and columns of $F$ as follows: for any two vertices $r_i$ and $r_j$, we write $r_i \preccurlyeq_r r_j$ if $r_i \rightarrow r_j$ in the semi-transitive orientation; and for any two vertices $c_i$ and $c_j$, we write $c_i \preccurlyeq_c c_j$ if $c_i \rightarrow c_j$ in the semi-transitive orientation. By \textit{Theorem~\ref{circR}}, it follows that $F$ does not have the circular-ones property. Hence, there exists at least one row-vertex, say $r_i$, that is not a circular interval of $\preccurlyeq_c$. Consequently, there exist four column-vertices $c_t \preccurlyeq_c c_u \preccurlyeq_c c_v \preccurlyeq_c c_w$ such that either $c_t, c_v \notin N_{CG(F)}(r_i)$ and $c_u, c_w \in N_{CG(F)}(r_i)$, or $c_t, c_v \in N_{CG(F)}(r_i)$ and $c_u, c_w \notin N_{CG(F)}(r_i)$. Since reversing all edges in a semi-transitive orientation also yields a semi-transitive orientation, there are four possible configurations, as illustrated in \textit{Figure~\ref{fp}}. In each case, the configuration induces either a shortcut or a directed cycle, which contradicts the assumption that the orientation is semi-transitive. This completes the proof of the lemma.
		\begin{figure}[h]
		\begin{minipage}{0.24\textwidth}
			\centering
			\begin{tikzpicture}[->,line cap=round,line join=round,>=triangle 45,x=1cm,y=1cm, scale=2]
				\begin{scriptsize}
					\node[draw, circle, fill = qqccqq] (1) at (1.5,2.6) {};
					\node[draw, circle, fill = qqzzcc] (5) at (1.8,2.4) {};
					\node[draw, circle, fill = qqzzcc] (c) at (1.8,2.75) {};
					\node[draw, circle, fill = qqzzcc] (d) at (1.8,3.15) {};
					\node[draw, circle, fill = qqzzcc] (b) at (1.8,2) {};
					
					\draw (c) -> (1);
					\draw (b) -> (1);
					\draw (5) -> (b);
					\draw (c) -> (5);
					\draw (d) -> (c);
				\end{scriptsize}
			\end{tikzpicture}
		\end{minipage}
		\hfill
		\begin{minipage}{0.24\textwidth}
			\centering
			\begin{tikzpicture}[->,line cap=round,line join=round,>=triangle 45,x=1cm,y=1cm, scale=2]
				\begin{scriptsize}
					\node[draw, circle, fill = qqccqq] (1) at (1.5,2.6) {};
					\node[draw, circle, fill = qqzzcc] (5) at (1.8,2.4) {};
					\node[draw, circle, fill = qqzzcc] (c) at (1.8,2.75) {};
					\node[draw, circle, fill = qqzzcc] (d) at (1.8,3.15) {};
					\node[draw, circle, fill = qqzzcc] (b) at (1.8,2) {};
					
					\draw (1) -> (c);
					\draw (b) -> (1);
					\draw (5) -> (b);
					\draw (c) -> (5);
					\draw (d) -> (c);
				\end{scriptsize}
			\end{tikzpicture}
		\end{minipage}
		\hfill
		\begin{minipage}{0.24\textwidth}
			\centering
			\begin{tikzpicture}[->,line cap=round,line join=round,>=triangle 45,x=1cm,y=1cm, scale=2]
				\begin{scriptsize}
					\node[draw, circle, fill = qqccqq] (1) at (1.5,2.6) {};
					\node[draw, circle, fill = qqzzcc] (5) at (1.8,2.4) {};
					\node[draw, circle, fill = qqzzcc] (c) at (1.8,2.75) {};
					\node[draw, circle, fill = qqzzcc] (d) at (1.8,3.15) {};
					\node[draw, circle, fill = qqzzcc] (b) at (1.8,2) {};
					
					\draw (1) -> (c);
					\draw (1) -> (b);
					\draw (5) -> (b);
					\draw (c) -> (5);
					\draw (d) -> (c);
				\end{scriptsize}
			\end{tikzpicture}
		\end{minipage}
			\begin{minipage}{0.24\textwidth}
				\centering
				\begin{tikzpicture}[->,line cap=round,line join=round,>=triangle 45,x=1cm,y=1cm, scale=2]
					\begin{scriptsize}
						\node[draw, circle, fill = qqccqq] (1) at (1.5,2.6) {};
						\node[draw, circle, fill = qqzzcc] (5) at (1.8,2.4) {};
						\node[draw, circle, fill = qqzzcc] (c) at (1.8,2.75) {};
						\node[draw, circle, fill = qqzzcc] (d) at (1.8,3.15) {};
						\node[draw, circle, fill = qqzzcc] (b) at (1.8,2) {};
						
						\draw (c) -> (1);
						\draw (1) -> (b);
						\draw (5) -> (b);
						\draw (c) -> (5);
						\draw (d) -> (c);
					\end{scriptsize}
				\end{tikzpicture}
			\end{minipage}
		\vspace{0.25cm}\caption{Shortcut possibilities}
		\label{fp}
	\end{figure}
	\end{proof}
	\begin{theorem}\label{scco}
		Let $G$ be a co-bipartite graph. $G$ is semi-transitive if and only if $M(G)$ has circularly compatible ones property.
	\end{theorem}
	\begin{proof}
		Let $G$ be a co-bipartite graph, and let $(X, Y)$ be the bipartition of the vertex set $V(G)$. Suppose that $G$ is semi-transitive. Assume, for the sake of contradiction, that $M(G)$ does not have the circularly compatible ones property. Then, by Theorem~\ref{cco}, $M(G)$ contains some matrix in $\FCCO^{\infty}$ as a configuration. However, by Lemma~\ref{cgf} and the hereditary nature of semi-transitive graphs, it follows that $G$ is not semi-transitive, a contradiction.
		
		Now suppose that $M(G)$ has the circularly compatible ones property. The idea of the proof is to orient the graph according to the circularly compatible ones biorder and to show that the resulting orientation is semi-transitive. Observe that the matrix $M(G)$ falls into one of the following cases:
		\begin{enumerate}
			\item $M(G)$ has either no trivial rows or no trivial columns.
			\item $M(G)$ has an all-$0$s row and an all-$0$s column.
			\item $M(G)$ has an all-$1$s row and an all-$1$s column.
		\end{enumerate}
		
		\textbf{Case 1:} Assume that $M(G)$ has either no trivial rows or no trivial columns. Without loss of generality, we may assume that $M(G)$ has no trivial rows (since $CG(M) \cong CG(M\trans)$). Therefore, by Theorems~\ref{cco} and~\ref{mco}, $M(G)$ admits a monotone circular biorder $(\preccurlyeq_r, \preccurlyeq_c)$. Let $\{r_1, r_2, \ldots, r_m\}$ and $\{c_1, c_2, \ldots, c_n\}$ denote the sets of rows and columns of $M$, listed in ascending order with respect to $\preccurlyeq_r$ and $\preccurlyeq_c$, respectively. By definition, we then have $X = \{r_1, r_2, \ldots, r_m\}$ and $Y = \{c_1, c_2, \ldots, c_n\}$. We define the orientation $\prec$ of the graph induced by the biorder $(\preccurlyeq_r, \preccurlyeq_c)$ as follows:
		\begin{itemize}
			\item For all $r_i, r_j \in X$, $r_i \prec r_j$ if and only if $r_i \preccurlyeq_r r_j$.
			\item For all $c_i, c_j \in Y$, $c_i \prec c_j$ if and only if $c_i \preccurlyeq_c c_j$.
			\item If $r_i = [d_i, e_i]_{\preccurlyeq_c}$ with $d_i \preccurlyeq_c e_i$, then $r_i \prec x$ for all $x$ such that $d_i \preccurlyeq_c x \preccurlyeq_c e_i$.
			\item If $r_i = [d_i, e_i]_{\preccurlyeq_c}$ with $e_i \prec_c d_i$, then $x \prec r_i$ for all $x \preccurlyeq_c e_i$ and $r_i \prec x$ for all $x$ with $d_i \preccurlyeq_c x$.
		\end{itemize}
		
		Before proceeding with the proof, we introduce some definitions. Let $r_i$ be a vertex that belongs to $X$. We call $r_i$ a \textit{linear-row vertex} if $r_i = [d_i, e_i]_{\preccurlyeq_c}$ with $d_i \preccurlyeq_c e_i$, and a \textit{circular-row vertex} if $r_i = [d_i, e_i]_{\preccurlyeq_c}$ with $e_i \prec_c d_i$.
		
		\textit{Claim~1.}\label{c1} \textit{If $r_i$ is a linear-row vertex and $r_j$ is a circular-row vertex, then $r_i \prec r_j$. Moreover, if $r_i = [d_i, e_i]_{\preccurlyeq_c}$ and $r_j = [d_j, e_j]_{\preccurlyeq_c}$, then $e_j \preccurlyeq_c e_i$.}
		
		Suppose, for the sake of contradiction, that $r_j \prec r_i$. This implies $r_j \prec_r r_i$. Let $r_i = [d_i, f_i]_{\preccurlyeq_c}^+$ and $r_j = [d_j, f_j]_{\preccurlyeq_c}^+$ be the unwrapped circular intervals, where $f_i = e_i$ and $f_j = e_j^+$ (see \textit{Definition}~\ref{ui}). By the definition of the monotone circular property, we have $f_j \preccurlyeq_c^+ f_i$, which implies $e_j^+ \preccurlyeq_c^+ e_i$, a contradiction. Therefore, $r_i \prec r_j$, as claimed.
		
		Hence, by the definition of the orientation, it follows that $r_q$ is a linear-row vertex for all $1 \le q \le i$, and a circular-row vertex for all $j \le q \le m$. This implies that $f_1 = e_1$ and $f_m = e_m^+$. Therefore, by the definition of the monotone circular property (see $(iii)$ of \textit{Definition~\ref{mc}}), it follows that $f_m \preccurlyeq_c^+ e_1^+$. Consequently, $e_m^+ \preccurlyeq_c^+ e_1^+$, which in turn implies that $e_m \preccurlyeq_c e_1$. Hence, we obtain the following ordering:
		\[
		e_j \preccurlyeq_c \cdots \preccurlyeq_c e_m \preccurlyeq_c e_1 \preccurlyeq_c \cdots \preccurlyeq_c e_i.
		\]
		Therefore, $e_j \preccurlyeq_c e_i$. This completes the proof of the claim.
		
		\textit{Claim~2.} \label{c2} 
		\textit{Let $v_1 \rightarrow v_2 \rightarrow \cdots \rightarrow v_p$, denoted by $\vec{S}$, be a directed path with respect to the orientation $\prec$. If $\vec{S}$ does not contain any linear-row vertex, then $\vec{S}$ induces a transitive tournament.}
		
	    The claim follows directly if $\vec{S}$ does not contain any vertex from $X$. Therefore, suppose that $\vec{S}$ contains at least one circular-row vertex. Since $\vec{S}$ trivially induces a transitive tournament when $p < 3$, assume that $p \ge 3$. Consider any three vertices $v_a, v_b, v_c$ in $\vec{S}$ such that $v_a \rightarrow v_b \rightarrow v_c$ for some $1 \le a < b < c \le p$. Assume toward a contradiction that $v_a \not\rightarrow v_c$. Let $v_a$ be a circular-row vertex, where $v_a = [d_a, e_a]_{\preccurlyeq_c}$. This implies that $v_c \in Y$. If $v_b$ is also a circular-row vertex, where $v_b = [d_b, e_b]_{\preccurlyeq_c}$, then by the monotone circular property, we have $d_a \preccurlyeq_c d_b \preccurlyeq_c v_c$. This implies that $v_a \rightarrow v_c$, a contradiction. If $v_b \in Y$, then by the definition of the orientation, we again have $v_a \rightarrow v_c$, a contradiction.
		
		Suppose that $v_a \in Y$. This implies that $v_c$ is a circular-row vertex, where $v_c = [d_c, e_c]_{\preccurlyeq_c}$. If $v_b$ is also a circular-row vertex, where $v_b = [d_b, e_b]_{\preccurlyeq_c}$, then by the monotone circular property, we have $v_a \preccurlyeq_c e_b \preccurlyeq_c e_c$. This implies that $v_a \rightarrow v_c$, a contradiction. If $v_b \in Y$, then by the definition of the orientation, we again have $v_a \rightarrow v_c$, a contradiction. Therefore, for any three vertices $v_a, v_b, v_c$ in $\vec{S}$ such that $v_a \rightarrow v_b \rightarrow v_c$ for some $1 \le a < b < c \le p$, we have $v_a \rightarrow v_c$. Hence, $v_1 \rightarrow v_2 \rightarrow \cdots \rightarrow v_p$ induces a transitive tournament. This completes the proof of the claim.
		
		We claim that the given orientation $\prec$ is semi-transitive. Suppose, for the sake of contradiction, that the given orientation is not semi-transitive. Let $\vec{S}$ be a shortcut path of minimal length, $\vec{S} = v_1 \rightarrow v_2 \rightarrow \cdots \rightarrow v_p$, where $v_1 \rightarrow v_p$ is the shortcut edge and $p \ge 4$. By \textit{Claim~2}, it is evident that $\vec{S}$ contains at least one linear-row vertex.
		
		Suppose that $k_1$ is the greatest integer in $[p]$ such that $v_{k_1}$ is a linear-row vertex. Then, by \textit{Claim~1} and by the definition of the orientation, $v_i$ is also a linear-row vertex for all $1 \le i \le k_1$. Similarly, let $k_2$ be the smallest integer in $[p]$ such that $v_{k_2} \in Y$. Therefore, by \textit{Claim~1} and the definition of the orientation, $v_j$ is a circular-row vertex for all $k_1 < j < k_2$.
		
		Now, suppose to the contrary that $k_2 > 3$. We claim that the path
		$v_1 \rightarrow v_{k_2 - 1} \rightarrow \cdots \rightarrow v_p$, denoted by $\vec{S'}$, induces a shortcut path with $v_1 \rightarrow v_p$ being the shortcut edge. Assume that $\vec{S'}$ contains more than three vertices. Then, for the sake of contradiction, suppose that $\vec{S'}$ induces a transitive tournament. Let $v_t$ be a vertex such that $k_2 \le t \le p$. If $v_t \in X$, then $v_l \rightarrow v_t$ for all $1 \le l < k_2$. Suppose $v_t \in Y$. Since $v_1 \rightarrow v_t$ and $v_{k_2 - 1} \rightarrow v_t$ for all $k_2 \le t \le p$, by the monotone circular property, we obtain the following ordering:
		\[
		d_1 \preccurlyeq_c \cdots \preccurlyeq_c d_{k_2 - 1} \preccurlyeq_c v_t \preccurlyeq_c e_1 \preccurlyeq_c \cdots \preccurlyeq_c e_{k_1}.
		\]
		
		By the definition of the orientation, $v_l \rightarrow v_t$ for all $1 \le l < k_2$ and $k_2 \le t \le p$. This implies that $\vec{S}$ induces a transitive tournament, a contradiction. Suppose that $\vec{S'}$ contains at most three vertices. Then $v_{k_2} = v_p$. Since $p = k_2 \ge 4$, and $v_1 \rightarrow v_p$ as well as $v_{k_2 - 1} \rightarrow v_p$, by the monotone circular property, we have $v_l \rightarrow v_p$ for all $1 \le l < k_2$. This again implies that $\vec{S}$ induces a transitive tournament, a contradiction. Therefore, $\vec{S'}$ must induce a shortcut path. However, since $k_2 > 3$, $\vec{S'}$ contains fewer vertices than $\vec{S}$, contradicting the choice of $\vec{S}$. Hence, $k_2 \le 3$.
		
		Now suppose that $k_3$ is the greatest integer in $[p]$ such that $v_{k_3} \in Y$. Assume, for the sake of contradiction, that $k_3 < p$. If $k_3 = k_2$, then by the definition of the monotone circular property and the second part of \textit{Claim~1}, we obtain the following ordering:
		\[
		d_1 \preccurlyeq_c \cdots \preccurlyeq_c d_{k_2 - 1} \preccurlyeq_c v_{k_2} \preccurlyeq_c e_{k_1 + 1} \preccurlyeq_c \cdots \preccurlyeq_c e_p \preccurlyeq_c e_1 \preccurlyeq_c \cdots \preccurlyeq_c e_{k_1}.
		\]
		This implies that $\vec{S}$ induces a transitive tournament, which is a contradiction.
		
		Now assume that $k_3 > k_2$. Then, by \textit{Claim~1}, it follows that $v_l$ is a circular-row vertex for all $k_3 < l \le p$. Therefore, by the definition of the monotone circular property and the second part of \textit{Claim~1}, we have $d_1 \preccurlyeq_c v_{k_3} \preccurlyeq_c e_{1}$. Hence, $v_1 \rightarrow v_{k_3}$. We claim that the path $v_1 \rightarrow v_2 \rightarrow \cdots \rightarrow v_{k_3}$, denoted by $\vec{S'}$, induces a shortcut path with $v_1 \rightarrow v_{k_3}$ being the shortcut edge. Assume that $k_3 \ge 4$. Suppose, for the sake of contradiction, that $\vec{S'}$ induces a transitive tournament. Let $v_t$ be a vertex for some $1 < t \le k_3$. If $v_t \in X$, then $v_t \rightarrow v_l$ for all $k_3 < l \le p$. Suppose $v_t \in Y$. Since $v_{k_3} \rightarrow v_{k_3 + 1}$, by the definition of the orientation and the monotone circular property, we obtain $v_t \rightarrow v_l$ for all $k_3 < l \le p$. This implies that $\vec{S}$ induces a transitive tournament, a contradiction.
		
		Now suppose that $k_3 \le 3$. This implies that $v_1$ is a linear-row vertex, $v_2, v_3 \in Y$, and $v_l$ is a circular-row vertex for all $4 \le l \le p$. By the definition of the monotone circular property and the second part of \textit{Claim~1}, it follows that $\vec{S}$ induces a transitive tournament, which is a contradiction. Therefore, $\vec{S'}$ induces a shortcut path. However, since $k_3 < p$, $\vec{S'}$ contains fewer vertices than $\vec{S}$, contradicting the choice of $\vec{S}$. Hence, $v_p \in Y$.
		
		As stated in \textit{Claim~2}, the path $v_3 \rightarrow v_4 \rightarrow \cdots \rightarrow v_p$ induces a transitive tournament, since it does not contain any linear-row vertex. Let $v_t$ be a vertex for some $3 \le t \le p$. If $v_2 \in Y$ or if $v_2$ is a circular-row vertex, then by the same argument, it follows that $v_2 \rightarrow v_t$. Since $v_1$ is a linear-row vertex, if $v_t \in X$, then $v_1 \rightarrow v_t$. Suppose $v_t \in Y$. Since $v_1 \rightarrow v_2$ and $v_1 \rightarrow v_p$, by the definition of the monotone circular property, we have $d_1 \preccurlyeq_c v_t \preccurlyeq_c e_1$. This implies that $v_1 \rightarrow v_t$, which in turn implies that $\vec{S}$ induces a transitive tournament, a contradiction.
		
		Now suppose that $v_2$ is a linear-row vertex. Since $v_1 \rightarrow v_p$ and $v_2 \rightarrow v_3$, by the definition of the monotone circular property, we have $d_1 \preccurlyeq_c d_2 \preccurlyeq_c v_t \preccurlyeq_c e_1 \preccurlyeq_c e_2$. Therefore, $v_1 \rightarrow v_t$ and $v_2 \rightarrow v_t$. This implies that $\vec{S}$ induces a transitive tournament, a contradiction. Since we have exhausted all possible cases, this contradicts the assumption that the orientation $\prec$ is not semi-transitive. Hence, the given orientation $\prec$ is semi-transitive. This completes the proof of this case.
		
		\textbf{Case 2: }Assume that $M(G)$ has an all-$0$s row and an all-$0$s column. Then, by the definition of the circularly compatible ones property, it follows that $M(G)$ has the circular-ones property for both rows and columns. Moreover, since it has an all-$0$s row and an all-$0$s column, by \textit{Theorem~\ref{c1cp}}, it follows that $M(G)$ has the consecutive-ones property for both rows and columns. Let $(\preccurlyeq_r, \preccurlyeq_c)$ be a circularly compatible-ones biorder such that $\{r_1, r_2, \ldots, r_m\}$ is an ascending order of $\preccurlyeq_r$, with $r_m$ being the all-$0$s row, and $\{c_1, c_2, \ldots, c_n\}$ is an ascending order of $\preccurlyeq_c$, with $c_n$ being the all-$0$s column. We define the orientation $\prec$ of the graph induced by the biorder $(\preccurlyeq_r, \preccurlyeq_c)$ as described in \textit{Case~1}.
		
		We claim that the given orientation $\prec$ is semi-transitive. Suppose to the contrary that the given orientation is not semi-transitive. Let $\vec{S}$ be a shortcut path of minimal length, $\vec{S} = v_1 \rightarrow v_2 \rightarrow \cdots \rightarrow v_p$, where $v_1 \rightarrow v_p$ is the shortcut edge and $p \ge 4$. If $\vec{S}$ does not contain either $r_m$ or $c_n$, then by \textit{Case~1}, it follows that $\vec{S}$ induces a transitive tournament. Therefore, assume that $\vec{S}$ contains either $r_m$ or $c_n$. By the definition of the orientation, it is evident that $\vec{S}$ cannot contain both $r_m$ and $c_n$, as there are no outgoing edges from either of these vertices. Moreover, whichever of the two $\vec{S}$ contains must be its sink vertex. 
		
		Suppose $v_p = c_n$. Then, by the definition of the orientation, it follows that $v_1 \in Y$. As there are no outgoing edges from the partition $Y$, it is evident that $\vec{S}$ induces a transitive tournament, a contradiction. Now suppose $v_p = r_m$. Then, by the definition of the orientation, we have $v_1 \in X$. Assume that there exists a vertex $v_i \in Y$ for some $2 \le i < p$. Since there are no outgoing edges from vertices in $Y$, it follows that every vertex $v_j$ with $i \le j \le p$ also belongs to $Y$. However, this contradicts the fact that $v_p \in X$. Therefore, $\vec{S}$ does not contain any vertex from the partition $Y$. Consequently, $\vec{S}$ induces a transitive tournament, a contradiction. Hence, the orientation $\prec$ is semi-transitive.
		
		\textbf{Case 3: }Assume that $M(G)$ has an all-$1$s row and an all-$1$s column. This implies that $G$ contains a universal vertex, say $v$. Therefore, by \textit{Theorem~\ref{neighbourhood}}, $G$ is semi-transitive if and only if $G - v$ is a permutation graph. Assume toward a contradiction that $G - v$ is not a permutation graph. Then, by \textit{Corollary~\ref{fisccirc}}, it follows that $G - v$ contains some graph $H \in \gc$ as an induced subgraph. Let $H'$ denote the subgraph of $G$ induced by $V(H) \cup \{v\}$.  
		
		If $H \cong \overline{C_{2k}}$ for some $k \ge 3$, then it can be easily verified that $M(H')$ contains the same configuration as $\MIast{k}$. If $H \cong G_1$, then $M(H')$ contains the same configuration as $\CoZtwo$. If $H \cong G_2$, then $M(H')$ contains the same configuration as $\Ztthree$. Finally, if $H \cong G_3$, then $M(H')$ contains $\MIast{3}$ as a configuration.  
		
		Therefore, from all the above cases, we conclude that $M(G)$ contains some matrix in $\FCCO$ as a configuration, which is a contradiction. This completes the proof of the theorem.
	\end{proof}
	
	\begin{figure}[h]
		
		\begin{minipage}{0.32\textwidth}
			\centering
			\begin{tikzpicture}[line cap=round,line join=round,>=triangle 45,x=1cm,y=1cm, scale=2]
				\begin{scriptsize}
					\node[draw, circle, fill = qqzzcc] (k-1) at (1.2,3.25) {};
					\draw[color=qqzzcc] (1,3.4) node {k-2};
					\node[draw, circle, fill = qqzzcc] (0) at (1.5,2.75) {};
					\node[draw, circle, fill = qqccqq] (3) at (1.8,3.25) {};
					\draw[color=qqccqq] (1.88,3.4) node {3};
					\node[draw, circle, fill = qqzzcc] (2) at (2.1,2.75) {};
					\draw[color=qqzzcc] (2.25,2.72) node {2};
					\node[draw, circle, fill = qqccqq] (1) at (1.8,2.25) {};
					\draw[color=qqccqq] (1.88,2.1) node {1};
					\node[draw, circle, fill = qqccqq] (k) at (0.9,2.75) {};
					\draw[color=qqccqq] (0.65,2.73) node {k - 1};
					\draw[color=black] (1.5,3.4) node {...};
					\node[draw, circle, fill = qqzzcc] (k+1) at (1.2,2.25) {};
					\draw[color=qqzzcc] (1.1,2.1) node {k};
					
					\draw (1) -- (2);
					\draw (3) -- (2);
					\draw (k-1) -- (k);
					\draw (k) -- (k+1);
					\draw (1) -- (k+1);
					\draw (1) -- (3);
					\draw (1) -- (k);
					\draw (k) -- (3);
					\draw (k+1) -- (2);
					\draw (k-1) -- (2);
					\draw (k-1) -- (k+1);
					\draw (k+1) -- (0);
					\draw (0) -- (2);
					\draw (k-1) -- (0);
				\end{scriptsize}
			\end{tikzpicture}\\{$CG(\MIast{k})$ for $k \geq 3$}
		\end{minipage}
		\begin{minipage}{0.32\textwidth}
			\centering
			\begin{tikzpicture}[line cap=round,line join=round,>=triangle 45,x=1cm,y=1cm, scale=2]
				\begin{scriptsize}
					\node[draw, circle, fill = qqzzcc] (k-1) at (1.2,3.25) {};
					\draw[color=qqzzcc] (1,3.4) node {k-2};
					\node[draw, circle, fill = qqzzcc] (0) at (1.5,2) {};
					\node[draw, circle, fill = qqccqq] (3) at (1.8,3.25) {};
					\draw[color=qqccqq] (1.88,3.4) node {3};
					\node[draw, circle, fill = qqzzcc] (2) at (2.1,2.75) {};
					\draw[color=qqzzcc] (2.25,2.72) node {2};
					\node[draw, circle, fill = qqccqq] (1) at (1.8,2.25) {};
					\draw[color=qqccqq] (1.88,2.1) node {1};
					\node[draw, circle, fill = qqccqq] (k) at (0.9,2.75) {};
					\draw[color=qqccqq] (0.65,2.73) node {k - 1};
					\draw[color=black] (1.5,3.4) node {...};
					\node[draw, circle, fill = qqzzcc] (k+1) at (1.2,2.25) {};
					\draw[color=qqzzcc] (1.1,2.1) node {k};
					
					\draw (k-1) -- (3);
					\draw (k-1) -- (0);
					\draw (k) -- (0);
					\draw (1) -- (3);
					\draw (0) -- (3);
					\draw (0) -- (1);
					\draw (0) -- (2);
					\draw (k) -- (3);
					\draw (k+1) -- (0);
					\draw (k+1) -- (3);
					\draw (k-1) -- (1);
					\draw (k) -- (1);
					\draw (k-1) -- (2);
					\draw (k+1) -- (2);
					\draw (k) -- (2);
					\draw (k-1) -- (k+1);
				\end{scriptsize}
			\end{tikzpicture}\\{$CG(\overline{\MIast{k}})$ for $k \geq 3$}
		\end{minipage}
		\begin{minipage}{0.32\textwidth}
			\centering
			\begin{tikzpicture}[line cap=round,line join=round,>=triangle 45,x=1cm,y=1cm, scale=2]
				\begin{scriptsize}
					\node[draw, circle, fill = qqzzcc] (6) at (1.05,2.15) {};
				\draw[color=qqzzcc] (0.85,2.12) node {6};
				\node[draw, circle, fill = qqzzcc] (5) at (1.5,1.9) {};
				\node[draw, circle, fill = qqzzcc] (8) at (1.5,2.5) {};
				\draw[color=qqzzcc] (1.5,2.63) node {8};
				\node[draw, circle, fill = qqzzcc] (7) at (1.95,2.15) {};
				\draw[color=qqzzcc] (2.1,2.12) node {7};
				\node[draw, circle, fill = qqccqq] (4) at (0.85,2.75) {};
				\draw[color=qqzzcc] (0.675,2.75) node {4};
				\node[draw, circle, fill = qqccqq] (2) at (1.95,3.35) {};
				\draw[color=qqccqq] (2.1,3.33) node {2};
				\node[draw, circle, fill = qqccqq] (3) at (2.15,2.75) {};
				\draw[color=qqccqq] (2.3,2.75) node {3};
				\node[draw, circle, fill = qqccqq] (1) at (1.05,3.35) {};
				\draw[color=qqccqq] (0.85,3.33) node {1};
				\draw[color=qqzzcc] (1.5,1.75) node {5};
					\draw (4) -- (1);
					\draw (1) -- (2);
					\draw (1) -- (3);
					\draw (1) -- (5);
					\draw (2) -- (3);
					\draw (2) -- (4);
					\draw (2) -- (5);
					\draw (2) -- (6);
					\draw (3) -- (5);
					\draw (3) -- (6);
					\draw (3) -- (7);
					\draw (4) -- (3);
					\draw (4) -- (6);
					\draw (5) -- (6);
					\draw (5) -- (7);
					\draw (5) -- (8);
					\draw (6) -- (7);
					\draw (6) -- (8);
					\draw (8) -- (7);
				
				\end{scriptsize}
			\end{tikzpicture}\\{$CG(\Ztwo)$}
		\end{minipage}\vspace{0.5cm}
		\begin{minipage}{0.32\textwidth}
			\centering
			\begin{tikzpicture}[line cap=round,line join=round,>=triangle 45,x=1cm,y=1cm, scale=2]
				\begin{scriptsize}
				 \node[draw, circle, fill = qqzzcc] (4) at (1.05,2.15) {};
				 \draw[color=qqzzcc] (0.85,2.12) node {7};
				 \node[draw, circle, fill = qqzzcc] (0) at (1.5,1.9) {};
				 \draw[color=qqzzcc] (1.5,1.75) node {8};
				 \node[draw, circle, fill = qqzzcc] (2) at (2.15,2.75) {};
				 \draw[color=qqzzcc] (2.3,2.75) node {5};
				 \node[draw, circle, fill = qqzzcc] (1) at (1.95,2.15) {};
				 \draw[color=qqzzcc] (2.1,2.12) node {6};
				 \node[draw, circle, fill = qqzzcc] (5) at (0.85,2.75) {};
				 \draw[color=qqzzcc] (0.675,2.75) node {4};
				 \node[draw, circle, fill = qqccqq] (c) at (1.95,3.35) {};
				 \draw[color=qqccqq] (2.1,3.33) node {2};
				 \node[draw, circle, fill = qqccqq] (d) at (1.5,3) {};
				 \draw[color=qqccqq] (1.5,3.15) node {3};
				 \node[draw, circle, fill = qqccqq] (b) at (1.05,3.35) {};
				 \draw[color=qqccqq] (0.85,3.33) node {1};
					
					\draw (4) -- (1);
					\draw (2) -- (1);
					\draw (5) -- (1);
					\draw (4) -- (2);
					\draw (5) -- (2);
					\draw (4) -- (5);
					\draw (c) -- (d);
					\draw (b) -- (d);
					\draw (b) -- (c);
					\draw (1) -- (c);
					\draw (2) -- (c);
					\draw (5) -- (b);
					\draw (1) -- (b);
					\draw (0) -- (1);
					\draw (0) -- (2);
					\draw (0) -- (4);
					\draw (0) -- (5);
					\draw (0) -- (b);
					\draw (0) -- (c);
					\draw (0) -- (d);
				\end{scriptsize}
			\end{tikzpicture}\\{$CG(\Ztthree)$}
		\end{minipage}
		\begin{minipage}{0.32\textwidth}
		\centering
		\begin{tikzpicture}[line cap=round,line join=round,>=triangle 45,x=1cm,y=1cm, scale=2]
			\begin{scriptsize}
				\node[draw, circle, fill = qqzzcc] (6) at (1.05,2.15) {};
			\draw[color=qqzzcc] (0.85,2.12) node {6};
			\node[draw, circle, fill = qqzzcc] (5) at (1.5,1.9) {};
			\node[draw, circle, fill = qqzzcc] (8) at (1.5,2.5) {};
			\draw[color=qqzzcc] (1.5,2.63) node {8};
			\node[draw, circle, fill = qqzzcc] (7) at (1.95,2.15) {};
			\draw[color=qqzzcc] (2.1,2.12) node {7};
			\node[draw, circle, fill = qqccqq] (4) at (0.85,2.75) {};
			\draw[color=qqzzcc] (0.675,2.75) node {4};
			\node[draw, circle, fill = qqccqq] (2) at (1.95,3.35) {};
			\draw[color=qqccqq] (2.1,3.33) node {2};
			\node[draw, circle, fill = qqccqq] (3) at (2.15,2.75) {};
			\draw[color=qqccqq] (2.3,2.75) node {3};
			\node[draw, circle, fill = qqccqq] (1) at (1.05,3.35) {};
			\draw[color=qqccqq] (0.85,3.33) node {1};
			\draw[color=qqzzcc] (1.5,1.75) node {5};
			\draw (4) -- (1);
			\draw (1) -- (2);
			\draw (1) -- (3);
			\draw (1) -- (5);
			\draw (1) -- (8);
			\draw (2) -- (3);
			\draw (2) -- (4);
			\draw (2) -- (5);
			\draw (2) -- (6);
			\draw (3) -- (5);
			\draw (3) -- (6);
			\draw (3) -- (7);
			\draw (4) -- (3);
			\draw (4) -- (6);
			\draw (5) -- (6);
			\draw (5) -- (7);
			\draw (5) -- (8);
			\draw (6) -- (7);
			\draw (6) -- (8);
			\draw (8) -- (7);
			\end{scriptsize}
		\end{tikzpicture}\\{$CG(\Zfive)$}
		\end{minipage}
		\begin{minipage}{0.32\textwidth}
			\centering
			\begin{tikzpicture}[line cap=round,line join=round,>=triangle 45,x=1cm,y=1cm, scale=2]
				\begin{scriptsize}
					\node[draw, circle, fill = qqzzcc] (4) at (1.05,2.15) {};
					\draw[color=qqzzcc] (0.85,2.12) node {7};
					\node[draw, circle, fill = qqzzcc] (0) at (1.5,1.9) {};
					\draw[color=qqzzcc] (1.5,1.75) node {8};
					\node[draw, circle, fill = qqzzcc] (2) at (2.15,2.75) {};
					\draw[color=qqzzcc] (2.3,2.75) node {5};
					\node[draw, circle, fill = qqzzcc] (1) at (1.95,2.15) {};
					\draw[color=qqzzcc] (2.1,2.12) node {6};
					\node[draw, circle, fill = qqzzcc] (5) at (0.85,2.75) {};
					\draw[color=qqzzcc] (0.675,2.75) node {4};
					\node[draw, circle, fill = qqccqq] (c) at (1.95,3.35) {};
					\draw[color=qqccqq] (2.1,3.33) node {2};
					\node[draw, circle, fill = qqccqq] (d) at (1.5,3) {};
					\draw[color=qqccqq] (1.5,3.15) node {3};
					\node[draw, circle, fill = qqccqq] (b) at (1.05,3.35) {};
					\draw[color=qqccqq] (0.85,3.33) node {1};
					
					\draw (4) -- (1);
					\draw (2) -- (1);
					\draw (5) -- (1);
					\draw (4) -- (2);
					\draw (5) -- (2);
					\draw (4) -- (5);
					\draw (5) -- (d);
					\draw (2) -- (d);
					\draw (2) -- (c);
					\draw (5) -- (b);
					\draw (1) -- (d);
					\draw (c) -- (d);
					\draw (b) -- (d);
					\draw (b) -- (c);
					\draw (0) -- (1);
					\draw (0) -- (2);
					\draw (0) -- (4);
					\draw (0) -- (5);
					\draw (0) -- (b);
					\draw (0) -- (c);
					\draw (0) -- (d);
				\end{scriptsize}
			\end{tikzpicture}\\{$CG(\CoZtwo)$}
		\end{minipage}
		\vspace{0.25cm}\caption{Minimal forbidden induced subgraphs for semi-transitive co-bipartite graphs, $\gs$. Note that the vertices with same colour induce a clique.}
		\label{gcob}
	\end{figure}
	\begin{lemma}\label{gco}
		If $F \in \FCCO^\infty$, then $CG(F)$ contains an induced subgraph isomorphic to some graph in $\gs$.
	\end{lemma}
	\begin{proof}
		Let $F$ be a matrix in $\FCCO^\infty$. If $F \in \{\MIast{k}, \overline{\MIast{k}} \colon k \ge 3\} \cup \{\Ztwo, \Ztthree, \CoZtwo, \Zfive\}$, then $CG(F)$ is isomorphic to a graph in $\gs$. If $F \in \{\Zfour, \CoZfour\}$, then $CG(F)$ is isomorphic to $CG(\Ztthree)$. This completes the proof of the lemma.
	\end{proof}
	\begin{lemma}\label{cnwr}
		All graphs in $\gs$ are minimal non-semi-transitive.
	\end{lemma}
	
	\begin{proof}
		Let $G$ be a graph in $\gs$. Observe that $M(G) \in \FCCO^\infty$. Therefore, by \textit{Theorems~\ref{cco}} and~\textit{\ref{scco}}, $G$ is a minimal non-semi-transitive graph.
	\end{proof}
	
	\begin{theorem}
		Let $G$ be a co-bipartite graph. Then $G$ is semi-transitive if and only if it contains none of the graphs in $\gs$ (depicted in Figure~\ref{gcob}) as an induced subgraph. Moreover, there exists a linear-time algorithm which, given a co-bipartite graph, determines whether it is semi-transitive.
	\end{theorem}
	
	\begin{proof}
		Assume that $G$ is semi-transitive. Suppose to the contrary that $G$ contains some graph in $\gs$ as an induced subgraph. Then, by \textit{Lemma~\ref{cnwr}} and the hereditary property of semi-transitive graphs, $G$ would be non-semi-transitive—a contradiction.
		
		Conversely, suppose that $G$ contains none of the graphs in $\gs$ as an induced subgraph. Assume, for the sake of contradiction, that $G$ is not semi-transitive. Then, by \textit{Theorem~\ref{scco}}, the matrix $M(G)$ does not have the circularly compatible ones property. Moreover, by \textit{Theorem~\ref{cco}}, $M(G)$ contains a matrix $F \in \FCCO^\infty$ as a configuration. Hence, by \textit{Lemma~\ref{gco}}, $G$ contains an induced subgraph isomorphic to a graph in $\gs$, a contradiction.
		
		Finally, note that for any co-bipartite graph $G$, $M(G)$ can be obtained in linear time. Hence, by \textit{Theorems~\ref{lta}} and~\textit{\ref{scco}}, the algorithmic part of the statement follows. This completes the proof.
	\end{proof}
The theorem provides a complete forbidden induced subgraph characterization of semi-transitive co-bipartite graphs. Specifically, it shows that a co-bipartite graph is semi-transitive if and only if it avoids a fixed finite family $\gs$ of induced subgraphs. Furthermore, this characterization leads to a linear-time algorithm for recognizing semi-transitive co-bipartite graphs, based on the circularly compatible ones property of the corresponding matrix representation.

\section{Conclusion and future work}\label{s5}

In this paper, we first establish a fundamental result: a co-bipartite graph is a circle graph if and only if it is a permutation graph. This equivalence leads to a minimal forbidden induced subgraph characterization of co-bipartite circle graphs. The principal contribution of this work is the explicit connection between the semi-transitivity of a co-bipartite graph and the circularly compatible ones property of its associated matrix. Using Safe’s characterization of matrices that satisfy the circularly compatible ones property, we derive a complete description of the minimal forbidden induced subgraphs for semi-transitive co-bipartite graphs, denoted by $\mathcal{G}_{co}$. Furthermore, this matrix-theoretic approach yields a linear-time recognition algorithm, providing an efficient solution to the recognition problem for this class of graphs. Several promising directions for future research emerge from this study:

\begin{enumerate}
	\item\textbf{Generalization to Multipartite Graphs:}
	The present work focuses on co-bipartite graphs, that is, graphs whose vertex set can be partitioned into two cliques. A natural direction for further research is the study of word-representability in graphs whose vertices can be partitioned into $k$ cliques for $k \ge 3$. An interesting open problem is whether the matrix-theoretic framework introduced in this work—most notably the circularly compatible ones property or appropriate generalizations of it—can be extended to provide structural and algorithmic characterizations for these broader classes of graphs.
	
	\item \textbf{Enumerative Combinatorics:}
	The explicit characterization of minimal forbidden induced subgraphs obtained in this work opens the door to an enumerative study of word-representable co-bipartite graphs. Such an investigation may lead to new combinatorial insights and reveal connections with known integer sequences or established combinatorial structures.
	
	\item \textbf{Extension of the Matrix Framework:}
	The successful application of the circularly compatible ones property to the characterization of semi-transitivity suggests that other matrix properties may serve a similar purpose in the study of word-representability for broader graph classes, such as chordal or weakly chordal graphs.
\end{enumerate}

\end{document}